\documentclass[a4paper,10pt]{scrartcl}
\usepackage[english]{babel}
\usepackage[utf8]{inputenc}
\usepackage{enumerate}
\usepackage{geometry}
\usepackage{amsfonts,amsmath,amssymb,amsopn}
\usepackage{amsthm}
\usepackage{mathtools}
\usepackage{stmaryrd}
\usepackage{nicefrac}
\usepackage{exscale}

\usepackage[numbers,square]{natbib}
\usepackage{url} 
\usepackage[colorlinks=true,pdfpagelabels,unicode]{hyperref}
\hypersetup{linkcolor=blue, urlcolor=blue, citecolor=blue} 

\allowdisplaybreaks 

\theoremstyle{plain}

\newtheoremstyle{theo}
	{3pt} 
	{3pt} 
	{\itshape} 
	{} 
		{\bfseries} 
	{\\} 
	{ } 
	{\thmname{#1}\thmnumber{ #2.}\thmnote{ - #3}} 
\theoremstyle{theo}

\newtheorem{definition}{Definition}[section]
\newtheorem{lemma}[definition]{Lemma}
\newtheorem{theorem}[definition]{Theorem}
\newtheorem{corollary}[definition]{Corollary}
\newtheorem{remark}[definition]{Remark}

\newenvironment{bew}{\begin{proof}[\bfseries Proof:]}{\end{proof}}

\DeclareMathOperator{\bomega}{\overline{\Omega}}
\DeclareMathOperator{\romega}{\partial\Omega}

\DeclareMathOperator{\intd}{d\!}

\DeclareMathOperator{\dive}{\nabla\cdot}
\DeclareMathOperator{\wto}{\rightharpoonup}

\newcommand{\epsi}{\varepsilon}

\usepackage[usenames,dvipsnames,svgnames,table]{xcolor}

\usepackage{tikz}
\usetikzlibrary{arrows,fadings,positioning,patterns,shapes}
\usepackage{pgfplots}

\newcommand{\Tme}{T_{max,\;\!\epsi}}
\newcommand{\GNI}{Gagliardo--Nirenberg inequality}

\newcommand{\into}[1]{\int_0^{#1}\!}

\newcommand{\intoT}{\into{T}}

\newcommand{\intomega}{\int_{\Omega}\!} 

\newcommand{\intoTomega}{\intoT\!\intomega}
\newcommand{\intinfomega}{\int_0^\infty\!\!\intomega}

\newcommand{\Lo}[1][1]{L^{#1}(\Omega)} 
\newcommand{\W}[1][1,2]{W^{#1}(\Omega)}

\newcommand{\LSp}[2]{L^{#1\;\!}\!\left(#2\right)} 
\newcommand{\LSploc}[2]{L_{loc}^{#1}\!\left(#2\right)} 
\newcommand{\WSp}[2]{W^{#1}\!\left(#2\right)}
\newcommand{\CSp}[2]{C^{#1}\!\left(#2\right)}

\newcommand{\R}{\mathbb{R}}
\newcommand{\N}{\mathbb{N}}

\newcommand{\nfrac}[2]{{\nicefrac{#1}{#2}}}

\author{Tobias Black\thanks{Institut f\"ur Mathematik, Universit\"at Paderborn, Warburger Str. 100, 33098 Paderborn, Germany; email: \mbox{tblack@math.upb.de}}}
\title{Global solvability of chemotaxis-fluid systems with nonlinear diffusion and matrix-valued sensitivities in three dimensions}
\date{}

\begin{document}
\maketitle
\begin{abstract}
\noindent
{\textbf{Abstract:} In this work we extend a recent result to chemotaxis fluid systems which include matrix-valued sensitivity functions $S(x,n,c):\Omega\times[0,\infty)^2\to\R^{3\times3}$ in addition to the porous medium type diffusion, which were discussed in the previous work. Namely, we will consider the system
\begin{align*}
\left\{
\begin{array}{r@{\,}c@{\,}c@{\ }l@{\quad}l@{\quad}l@{\,}c}
n_{t}&+&u\cdot\!\nabla n&=\Delta n^m-\nabla\!\cdot(nS(x,n,c)\nabla c),\ &x\in\Omega,& t>0,\\
c_{t}&+&u\cdot\!\nabla c&=\Delta c-c+n,\ &x\in\Omega,& t>0,\\
u_{t}&+&(u\cdot\nabla)u&=\Delta u+\nabla P+n\nabla\phi,\ &x\in\Omega,& t>0,\\
&&\nabla\cdot u&=0,\ &x\in\Omega,& t>0,
\end{array}\right.
\end{align*}
in a bounded domain $\Omega\subset\mathbb{R}^3$ with smooth boundary. Assuming that $m\geq1$, $\alpha\geq0$ satisfy $m+\alpha>\frac43$, that the matrix-valued function $S(x,n,c):\Omega\times[0,\infty)^2\to\R^{3\times3}$ satisfies $|S(x,n,c)|\leq\frac{S_0}{(1+n)^{\alpha}}$ for some $S_0>0$ and suitably regular nonnegative initial data, we show that the corresponding no-flux-Dirichlet boundary value problem emits at least one global very weak solution. Upon comparison with results for the fluid-free system this condition appears to be optimal.
Moreover, imposing a stronger condition for the exponents $m$ and $\alpha$, i.e. $m+2\alpha>\frac{5}{3}$, we will establish the existence of at least one global weak solution in the standard sense.
\noindent
}\\[0.1cm]


\end{abstract}


\newpage
\section{Introduction}\label{sec1:intro}
The coupling of chemotaxis, the biological phenomenon of directed movement of cells in response to a signal chemical present in the neighborhood of the organism, to the Navier-Stokes-fluid-equations, and thereby including interplay between cells, chemical and fluid surrounding, has been of increasing interest in the last decade. Studies on broadcast spawning indicate the influence this coupling can have on the migration process (\cite{coll1994chemical,miller1985demonstration}). Particular attention has thus been devoted to the question whether results known for the classical Keller--Segel-system (\cite{KS70},\cite{Ho03}) can be transferred to the setting incorporating this fluid interaction. A distinct feature of the Keller--Segel model (even without fluid) is its possibility to capture the emergence of patterns arising from the aggregation of bacteria, which on the solution level of the corresponding PDE system
\begin{align}\label{KS}
n_t=\nabla\cdot\big(D(n)\nabla n- S(n,c)\nabla c\big)&& c_t=\Delta c-c+n,
\end{align} 
with $n(x,t)$ denoting the cell density and $c(x,t)$ the signal concentration, can be observed as blow-up of solutions. Correspondingly, the significance of obtaining results proving or excluding the possibility of blow-up have been a very important concern of the literature. For an extensive overview of results we refer the reader to the survey \cite{BBWT15}. For the Keller--Segel system of the form in \eqref{KS} the quantity governing the behavior has been identified to be the growth ration of $\frac{S(n)}{D(n)}$, with its critical number given by $\frac{2}{N}$ and $N$ being the space dimension (see \cite{TaoWin-quasilinear_JDE12} and references therein). In fact, the sufficient conditions for blow-up to be excluded in the corresponding Neumann-boundary value problem in a smooth domain $\Omega\subset\R^N$ the classical solutions emerging from suitably regular initial data remain bounded for all times, whenever
\begin{align*}
\frac{S(n)}{D(n)}\leq C(n+1)^\beta\quad\text{for all }n\geq0\text{ with some }C>0\text{ and }\beta<\frac{2}{N}.
\end{align*}
On the other hand in \cite{Win-volume-filling-MMAS10} smooth solutions blowing-up in either infinite or finite time have been shown to exist under the assumption of 
\begin{align*}
\frac{S(n)}{D(n)}\geq Cn^\gamma\quad\text{for all }n>1\text{ with some }C>0\text{ and }\gamma>\frac{2}{N}.
\end{align*}
(Finite time blow-up has also been witnessed in \cite{CieslakStinner-JDE15}) Especially, considering cell diffusion as covered by variants of the porous medium operator, but nondegenerate, i.e. $D(n)\equiv m(n+1)^{m-1}$, and a sensitivity functions satisfying $S(n)\equiv(1+n)^{1-\alpha}$, the condition for finite time blow-up to be excluded in \eqref{KS} can be expressed as $m+\alpha>\frac{2N-2}{N}$, which will act as our comparison point for conditions arising in the setting with fluid. For the systems incorporating fluid interaction and signal production
\begin{align}\label{CTNS}
\left\{
\begin{array}{r@{\,}c@{\,}c@{\ }l@{\quad}l}
n_{t}&+&u\cdot\!\nabla n&=\nabla\cdot\big(D(n)\nabla n-nS(x,n,c)\nabla c\big),\\
c_{t}&+&u\cdot\!\nabla c&=\Delta c-c+n,\\
u_{t}&+&\kappa(u\cdot\nabla)u&=\Delta u+\nabla P+n\nabla\phi,\\
&&\nabla\cdot u&=0,
\end{array}\right.
\end{align}
where $S$ may be a tensor-valued function, $u$ now denotes the fluid-velocity field, $P$ the corresponding pressure and $\phi$ is a given gravitational potential, however, the literature is not as rich and mostly focuses either on the case $D(n)\equiv1$ or on $S(x,n,c)\equiv 1$. (A more common variant of \eqref{CTNS} is concerned with signal consumption and was proposed by \cite{tuval2005bacterial}. For this setting the results are a bit more extensive and an overview of known results in three-dimensional domains can be found in the references of \cite{TB2017_nonlindiff}.) Let us briefly recapitulate the recent developments for porous medium type diffusion $D(n)=m n^{m-1}$. In the case of $m=1$ (i.e. linear diffusion) and tensor-valued $S(x,n,c)$ satisfying $|S(x,n,c)|\leq(1+n)^{-\alpha}$ global weak solutions were shown to exist for $\alpha\geq\frac{3}{7}$ (\cite{LiuWang-GlobWeak-JDE17}) and global very weak solutions were established whenever $\alpha>\frac{1}{3}$ (\cite{Wang17-globweak-ksns}). In space dimension $N=2$ the optimal condition $\alpha>0$ can even be reached with global bounded classical solutions (\cite{WangWinXiang-Pisa18}). If we simplify to Stokes-fluid ($\kappa=0$ in \eqref{CTNS}) instead of full Navier--Stokes-fluid, more regular solutions can also achieved in dimension $N=3$, as indicated by the recent studies on bounded classical solutions in \cite{Win18-alphastokesreg}. On the other hand, in the case of $S(x,n,c)\equiv 1$ (i.e. $\alpha=0$) and $m>1$ global weak solutions were obtained first for $m>2$ in \cite{ZhengJ-globweak-JDE17} and more recently for $m>\frac{5}{3}$ in \cite{TB2017_nonlindiff}, were also global very weak solutions were shown to exist whenever $m>\frac{4}{3}$. The results concerning $N=3$ and Navier--Stokes-fluid can be illustrated by the following picture.

\begin{center}
\begin{tikzpicture}
  \node[below] at (2,0){1};
  \node[below] at (4,0){2};
  \node[left] at (0,2){1};
  \node[left] at (0,4){2};
  \node[left] at (0,8/3){$\frac43$};
  \node[left] at (0,10/3){$\frac53$};
  \node[below] at (8/3,0){$\frac43$};
  \node[below] at (6/7,-0.1){$\frac37$};

\scope
\clip[postaction={fill=blue!30}] (6/7,2.04)--(5.3,2.04)--(5.3,1.96)--(6/7,1.96)--cycle; 
\endscope
\scope
\clip[postaction={fill=blue!30}] (-0.04,10/3)--(0.04,10/3)--(0.04,5.3)--(-0.04,5.3)--cycle; 
\endscope
\scope
\clip[postaction={fill=red!40}] (-0.04,10/3)--(0.04,10/3)--(0.04,8/3)--(-0.04,8/3)--cycle; 
\endscope
\scope
\clip[postaction={fill=red!40}] (2/3-0.015,2.04)--(6/7,2.04)--(6/7,1.96)--(2/3+0.035,1.96)--cycle; 
\endscope  

  \draw[thin,gray,dashed] (6/7,2) -- (6/7,0);
  \draw[thin] (0,8/3) -- (8/3,0);
  \draw[thin,->] (-0.3,0) -- (5.35,0) coordinate[label = {below:$\alpha$}] (xmax);
  \draw[thin,->] (0,-0.3) -- (0,5.35) coordinate[label = {left:$m$}] (ymax);
  \draw[thin] (0,8/3) -- (8/3,0);
  \draw[gray,dashed,thin] (5.2,2) -- (2/3,2);

 \node[right,label={[shift={(2.0,-0.13)},anchor=east]}] at (5.7,5){Global existence of};
  \node[right,xscale=2,rectangle,fill=red!40] at (5.5,4.5){};
 \node[align=left,right] at (6,4.5){very weak solution}; 
  \node[right,xscale=2,rectangle,fill=blue!30,] at (5.5,4){};
 \node[align=left,right] at (6,4) {weak solution};
 \draw[thin,opacity=0.5,color=lightgray, step=2/3] (0,0) grid (5.3,5.3);
  \draw[thin] (2,-0.1) -- (2,0.1); 
  \draw[thin] (4,-0.1) -- (4,0.1);
  \draw[thin] (-0.1,2) -- (0.1,2);
  \draw[thin] (-0.1,4) -- (0.1,4);
  \draw[thin] (-0.05,8/3) -- (0.05,8/3);
  \draw[thin] (8/3,-0.05) -- (8/3,0.05);
  \draw[thin] (-0.05,10/3) -- (0.05,10/3);
  \draw[thin] (6/7,-0.05) -- (6/7,0.05);
\node[align=center, xshift=-5em,yshift=-0.8em] (title) 
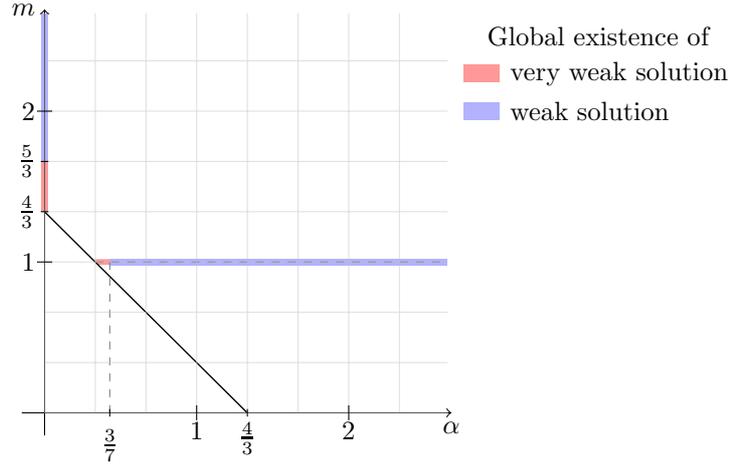
    at (current bounding box.south) {\footnotesize{Fig. 1.1: Overview of global existence with fluid interaction prior to this work}};  
\end{tikzpicture}
\end{center}

Comparing the either or cases above one expects global very weak solutions to exist for all $m\geq1$ and $\alpha\geq0$ satisfying $m+\alpha>\frac{4}{3}$. However, connecting the currently known limit cases for weak solutions in the standard sense to exist, leads to a line which appears to have a rather unnatural slope, posing the question whether the current condition $m=1$ and $\alpha>\frac{3}{7}$ is critical in $\alpha$ for global weak solutions to exist. Our main interest thereby consists in extracting a priori estimates from the sparse information provided by the system, which, most importantly, captures optimal conditions on $m\geq1$ and $\alpha\geq0$.

\noindent{\textbf{Main results.} 
Suppose that $\Omega\subset\R^3$ is a bounded domain with smooth boundary, $m\geq1$ and that for some $\alpha\geq0$ and $S_0>0$ the matrix-valued sensitivity function $S\in\CSp{2}{\bomega\times[0,\infty)^2;\R^{3\times3}}$ satisfies 
\begin{align}\label{Sprop}
|S(x,n,c)|\leq\frac{S_0}{(1+n)^{\alpha}}\quad\text{for all }x\in\bomega, n\geq0\text{ and }c\geq0.
\end{align}
Under these assumptions we will consider 
\begin{align}\label{CTndS}
\left\{
\begin{array}{r@{\,}c@{\,}c@{\ }l@{\quad}l@{\quad}l@{\,}c}
n_{t}&+&u\cdot\!\nabla n&=\Delta n^m-\nabla\!\cdot(nS(x,n,c)\nabla c),\ &x\in\Omega,& t>0,\\
c_{t}&+&u\cdot\!\nabla c&=\Delta c-c+n,\ &x\in\Omega,& t>0,\\
u_{t}&+&(u\cdot\nabla)u&=\Delta u+\nabla P+n\nabla\phi,\ &x\in\Omega,& t>0,\\
&&\nabla\cdot u&=0,\ &x\in\Omega,& t>0,
\end{array}\right.
\end{align}
complemented with boundary conditions
\begin{align}\label{BC}
\Big(\nabla n^m(x,t)-n(x,t)S\big(x,n(x,t),c(x,t)\big)\nabla c(x,t)\Big)\cdot\nu&=0,\nonumber\\\nabla c(x,t)\cdot\nu=0\quad\text{and}\quad u(x,t)&=0\qquad\text{for }x\in\romega\text{ and }t>0,
\end{align}
and initial conditions
\begin{align}\label{IC}
n(x,0)=n_0(x),\quad c(x,0)=c_0(x),\quad u(x,0)=u_0(x),\quad x\in\Omega,
\end{align}
where the gravitational potential $\phi$ is assumed to satisfy
\begin{align}\label{phidef}
\phi\in\W[2,\infty].
\end{align}
Prescribing initial data which satisfy the conditions
\begin{align}\label{IR}
\left\{\begin{array}{r@{\,}l}
n_0&\in\CSp{\gamma}{\bomega}\quad\text{for some }\gamma>0\quad\text{with } n_0\geq0\text{ in }\Omega\text{ and }n_0\not\equiv0,\\
c_0&\in\W[1,\infty]\quad\text{ with }c_0\geq0\text{ in }\bomega,\text{ and }c_0\not\equiv0,\\
u_0&\in\WSp{2,2}{\Omega;\R^3}\cap W_0^{1,2}\left(\Omega;\R^3\right)\quad\text{such that } \nabla\cdot u_0=0,
\end{array}\right.
\end{align}
we obtain the following main results.

\begin{theorem}\label{theo:1}
Let $\Omega\subset\R^3$ be a bounded domain with smooth boundary. Suppose that $m\geq1$ and $\alpha\geq0$ satisfy $m+2\alpha>\frac{5}{3}$. Moreover, assume $S\in\CSp{2}{\bomega\times[0,\infty)^2;\R^{3\times 3}}$ fulfills \eqref{Sprop} with some $S_0>0$ and that $n_0,c_0$ and $u_0$ comply with \eqref{IR}. Then \eqref{CTndS}--\,\eqref{IC} admits at least one global weak solution in the sense of Definition \ref{def:weak-solution} below.
\end{theorem}

\begin{remark}
For the linear diffusion case $m=1$ Theorem \ref{theo:1} provides the existence of a global weak solution for $\alpha>\tfrac{1}{3}$, extending the results of \cite{Wang17-globweak-ksns} and \cite{LiuWang-GlobWeak-JDE17}, which provided the existence of a global very weak solution for $\alpha>\tfrac13$ and a global weak solution for $\alpha>\tfrac37$, respectively.
\end{remark}

If we merely prescribe $m+2\alpha\leq\frac{5}{3}$, we have to weaken the solution concept in order to verify the existence of global solutions -- which is due to the obtainable a priori information being so weak that we have to consider a sublinear functional of $n$ for our testing methods.

\begin{theorem}\label{theo:2}
Let $\Omega\subset\R^3$ be a bounded domain with smooth boundary. Suppose that $m\geq1$ and $\alpha\geq0$ satisfy $m+\alpha>\frac{4}{3}$. Moreover, assume $S\in\CSp{2}{\bomega\times[0,\infty)^2;\R^{3\times 3}}$ fulfills \eqref{Sprop} and that $n_0,c_0$ and $u_0$ comply with \eqref{IR}. Then \eqref{CTndS}--\,\eqref{IC} admits at least one global very weak solution $(n,c,u)$ in the sense of Definition \ref{def:very_weak_sol} below. In particular, this global very weak solution satisfies
\begin{align*}
n&\in \LSploc{2(m+\alpha)-\frac{4}{3}}{\bomega\times[0,\infty)},\quad c\in \LSploc{2}{[0,\infty);\W[1,2]},\quad u\in \LSploc{2}{[0,\infty);W_{0,\sigma}^{1,2}(\Omega;\R^3)},
\end{align*}
and
\begin{align*}
\intomega n(\cdot,t)=\intomega n_0\quad\text{for a.e. }t>0.
\end{align*}
\end{theorem}

Since we are considering Navier--Stokes-fluid, smooth global solutions can not be expected. However, it could be expected that the very weak solutions obtained for $m+\alpha>\frac{4}{3}$ may in fact become smooth solutions after some waiting time. Effects of this kind have, in more generous setting featuring signal consumption instead of production, been observed in e.g. \cite{win15_chemonavstokesfinal}.

Illustrating the diagram of before once more with the new results, we obtain the following figure, which neatly fits together with the expectations we obtained from Figure 1.1.
\begin{center}
\begin{tikzpicture}
  \node[below] at (2,0){1};
  \node[below] at (4,0){2};
  \node[left] at (0,2){1};
  \node[left] at (0,4){2};
  \node[left] at (0,8/3){$\frac43$};
  \node[left] at (0,10/3){$\frac53$};
  \node[below] at (8/3,0){$\frac43$};
  \node[below] at (6/7,-0.1){$\frac37$};

  \scope
       \clip[postaction={fill=blue!30}] (5.3,2) -- (2/3,2) -- (0,10/3) -- (0,5.3) -- (5.3,5.3) --(5.3,2);
  \endscope
  \scope
 \clip[postaction={fill=red!40}] (0,8/3) -- (0,10/3) -- (2/3,2) -- (0,8/3);
  \endscope
 
 \draw[thin] (0,8/3) -- (8/3,0);
  \draw[thin] (2/3,2) -- (5.3,2);
   \draw[thin,gray,dashed] (6/7,2) -- (6/7,0);
    \draw[thin,->] (-0.3,0) -- (5.35,0) coordinate[label = {below:$\alpha$}] (xmax);
  \draw[thin,->] (0,-0.3) -- (0,5.35) coordinate[label = {left:$m$}] (ymax);
    \draw[thin] (0,8/3) -- (8/3,0);
  \draw[gray,dashed,thin] (0,10/3) -- (2/3,2);
  \draw[gray,dashed,thin] (5.2,2) -- (2/3,2);
 \draw[thin,opacity=0.5,color=lightgray, step=2/3] (0,0) grid (5.3,5.3);
  \draw[thin] (2,-0.1) -- (2,0.1);
  \draw[thin] (4,-0.1) -- (4,0.1);
  \draw[thin] (-0.1,2) -- (0.1,2);
  \draw[thin] (-0.1,4) -- (0.1,4);
  \draw[thin] (-0.05,8/3) -- (0.05,8/3);
  \draw[thin] (8/3,-0.05) -- (8/3,0.05);
  \draw[thin] (-0.05,10/3) -- (0.05,10/3);
  \draw[thin] (6/7,-0.05) -- (6/7,0.05);
 \node[right,label={[shift={(2.0,-0.13)},anchor=east]}] at (5.7,5){Global existence of};
 \node[right,xscale=2,rectangle,fill=gray!50] at (5.5,4.5){};
  \node[right,xscale=2,rectangle,fill=red!40] at (5.5,4.5){};
 \node[align=left,right] at (6,4.5){very weak solution}; 
 \node[right,xscale=2,rectangle,fill=black!60,] at (5.5,4){};
  \node[right,xscale=2,rectangle,fill=blue!30,] at (5.5,4){};
 \node[align=left,right] at (6,4) {weak solution};

\node[align=center, xshift=-5em, yshift=-0.8em] (title) 
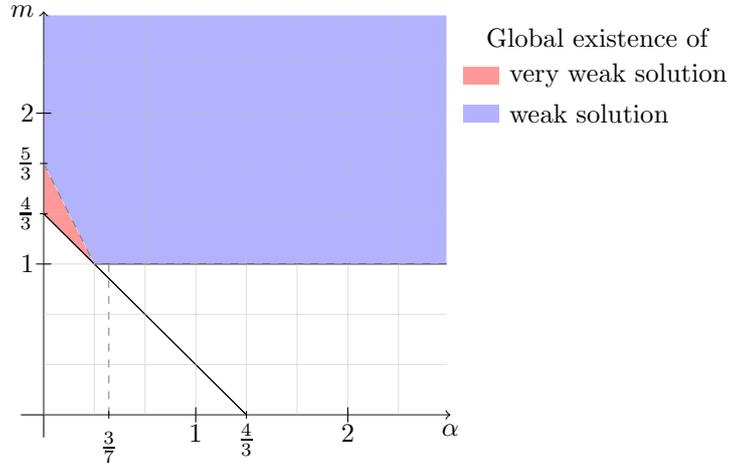
    at (current bounding box.south) {\footnotesize{Fig. 1.2: Overview of with Thm.1.1. and Thm.1.2.}};  
\end{tikzpicture}
\end{center}
\noindent{\textbf{Mathematical difficulties.} 
The absence of any energy-functional in this setting incorporating both fluid interaction and signal production, is one of the main difficulties in obtaining estimates optimal with respect to $m$ and $\alpha$. Most of the problems resulting from this lack of an energy estimate can be combated by utilizing similar methods as displayed in \cite{Wang17-globweak-ksns} and our previous work \cite{TB2017_nonlindiff}, but even greater care has to be taken when trying to derive information on gradient terms and combined quantities without tightening the scope for $m$ and $\alpha$. After regularizing the problem in a suitable fashion, a functional of the form
\begin{align*}
\intomega(n_\epsi+\epsi)^{m+2\alpha-1}(\cdot,t)+\intomega c_\epsi^2(\cdot,t),\quad t>0,
\end{align*}
(which for small values of $m$ and $\alpha$ is of sublinear growth with respect to $n$) will be the main cornerstone of our analysis and will also provide bounds on $\int_t^{t+1}\!\intomega\big|\nabla(n_\epsi+\epsi)^{m+\alpha-1}\big|^2$ as well as $\int_t^{t+1}\!\intomega|\nabla c_\epsi|^2$ (Lemma \ref{lem:bounds}), which by the \GNI\ can be refined into more spatio-temporal regularity information on $n_\epsi$ (Lemma \ref{lem:st-bound-nepsi}). Carefully combining these estimates with standard arguments for the Navier--Stokes-subsystem will enable us to conclude from compactness arguments the existnce of the desired limit object. (Lemma \ref{lem:convergences}). Depending on the size of $m$ and $\alpha$ the convergence properties can be relied on to conclude Theorems \ref{theo:1} and \ref{theo:2}.

\setcounter{equation}{0} 
\section{The notions of weak and very weak solutions}\label{sec2:soldef}
Let us start by laying out the except formulations of the different concepts of solvability we are going to discuss. The notion of very weak solvability present in Theorem \ref{theo:2} is adapted from the related works in \cite{win15_chemorot,Wang17-globweak-ksns,TB2017_nonlindiff} and the main distinguishing aspect when comparing to the standard notion of weak solvability is the fact that the first component of the system only has to satisfy a supersolution property. To be more precise, we require the following
\begin{definition}\label{def:weak_super_sol}
Let $\Phi\in\CSp{2}{[0,\infty)}$ be a nonnegative function satisfying $\Phi'>0$ on $(0,\infty)$. 
Assume that $n_0\in\Lo[\infty]$ is nonnegative and that $\Phi(n_0)\in\Lo[1]$. Moreover, let $S\in\CSp{2}{\bomega\times[0,\infty)^2;\R^{3\times3}}$ satisfy \eqref{Sprop} for some $S_0>0$ and $\alpha\geq0$. Suppose that $c\in\LSploc{2}{[0,\infty);\W[1,2]}$ and $u\in\LSploc{1}{[0,\infty); W_0^{1,1}\left(\Omega;\R^3\right)}$ with $\dive u\equiv 0$ in $\mathcal{D}'\big(\Omega\times(0,\infty)\big)$. 
The nonnegative measurable function $n:\Omega\times(0,\infty)\to\R$ satisfying $n\in\LSploc{1}{[0,\infty);\W[1,1]}$ will be named a global weak $\Phi$--supersolution of the initial-boundary value problem
\begin{align}\label{eq:system_very_weak_supersol}
\left\{
\begin{array}{r@{\,}c@{\,}r@{\ }l@{\quad}l@{\quad}l@{\,}c}
n_{t}&+&u\cdot\!\nabla n&=\Delta n^m-\nabla\!\cdot\big(nS(x,n,c)\nabla c\big),\ &x\in\Omega,& t>0,\\
&&\frac{\partial n}{\partial \nu}&=0,\ &x\in\romega,& t>0,\\
&&n(x,0)&=n_0(x),\ &x\in\Omega,&
\end{array}\right.
\end{align}
if 
\begin{align}\label{eq:very_weak_supersol_regularity}
\Phi(n),\text{ and }\ \Phi''(n)n^{m-1}|\nabla n|^2\text{ belong to }\LSploc{1}{\bomega\times[0,\infty)},&\nonumber\\
\Phi'(n)n^{m-1}\nabla n,\text{ and }\ \Phi(n)u\text{ belong to }\LSploc{1}{\bomega\times[0,\infty);\R^3},&\\
\Phi'(n)n\text{ belongs to }\LSploc{2}{\bomega\times[0,\infty)},\text{ and }\ \Phi''(n)n\nabla n\text{ belongs to }&\LSploc{2}{\bomega\times[0,\infty);\R^3},\nonumber
\end{align}
and if for each nonnegative $\varphi\in C_0^\infty\left(\bomega\times[0,\infty)\right)$ with $\frac{\partial\varphi}{\partial\nu}=0$ on $\romega\times(0,\infty)$, the inequality
\begin{align}\label{eq:very_weak_supersol}
-\intinfomega&\Phi(n)\varphi_t-\intomega\Phi(n_0)\varphi(\cdot,0)\nonumber\\
&\geq -\,m\intinfomega\Phi''(n)n^{m-1}|\nabla n|^2\varphi-m\intinfomega\Phi'(n)n^{m-1}(\nabla n\cdot\nabla \varphi)\\&\qquad\quad\,+\intinfomega\Phi''(n)n\big(\nabla n\cdot S(x,n,c)\nabla c\big)\varphi+\intinfomega\Phi'(n)n\big(S(x,n,c)\nabla c\cdot\nabla\varphi\big)\nonumber\\&\qquad\qquad\quad+\intinfomega\Phi(n)(u\cdot\nabla\varphi)\nonumber
\end{align}
is satisfied.
\end{definition}

Let us briefly remark on the test function we will use later on. For $m\geq1$ and $\alpha\geq0$ satisfying the conditions $m+\alpha>\frac{4}{3}$ and $m+2\alpha<2$ we will consider $\Phi(s)\equiv(s+1)^{m+2\alpha-1}$. Due to $m+2\alpha-1$ our main intention in the coming section will be to obtain a priori bounds which allow for the conclusion that $n^{m+2\alpha-1}\in\LSploc{2}{[0,\infty);\W[1,2]}$. Combining this with suitable regularity information on the other solution components is sufficient to determine that all of the integrals appearing in the supersolution property above are well defined (see als Corollary \ref{cor:1} and Lemma \ref{lem:sol-prop-n} below).

Complementing Definition \ref{def:weak_super_sol} with the standard properties of weak solvability for the remaining subproblems of \eqref{CTndS} will lead us to the following notion of global very weak solutions.

\begin{definition}\label{def:very_weak_sol}
A triple $(n,c,u)$ of functions
\begin{align*}
\begin{array}{r@{\ }c@{\ }l}
n&\in&\LSploc{1}{\bomega\times[0,\infty)},\\
c&\in&\LSploc{2}{[0,\infty);\W[1,2]},\\
u&\in& L_{loc}^1\big([0,\infty); W_{0}^{1,1}\!\left(\Omega;\R^3\right)\!\big),
\end{array}
\end{align*}
satisfying $n\geq0$ and $c\geq0$ in $\bomega\times[0,\infty)$, $cu\in\LSploc{1}{\bomega\times[0,\infty)}$, as well as $u\otimes u\in\LSploc{1}{\bomega\times[0,\infty);\R^{3\times3}}$ will be called a global very weak solution of \eqref{CTndS}--\,\eqref{IC}, if 
\begin{align*}
\intomega n(\cdot,t)\leq\intomega n_0\quad\text{for a.e. }t>0,
\end{align*}
if $\nabla\cdot u=0$ in $\mathcal{D}'\big(\Omega\times(0,\infty)\big)$, if the equality
\begin{align}\label{eq:very_weak_sol1}
-\intinfomega c\varphi_t-\intomega c_0\varphi(\cdot,0)=-\intinfomega\nabla c\cdot\nabla \varphi-\intinfomega c\varphi+\intinfomega n\varphi+\intinfomega c(u\cdot\nabla\varphi)
\end{align}
holds for all $\varphi\in\LSp{\infty}{\Omega\times(0,\infty)}\cap\LSp{2}{(0,\infty);\W[1,2]}$ with $\varphi_t\in\LSp{2}{\Omega\times(0,\infty)}$, which are compactly supported in $\bomega\times[0,\infty)$, if
\begin{align}\label{eq:very_weak_sol2}
-\intinfomega u\cdot\psi_t-\intomega u_0\cdot\psi(\cdot,0)=-\intinfomega\nabla u\cdot\nabla\psi+\intinfomega(u\otimes u)\cdot\nabla\psi+\intinfomega n\nabla\phi\cdot\psi
\end{align}
is fulfilled for all $\psi\in C_0^\infty\left(\Omega\times[0,\infty);\R^3\right)$ with $\nabla\cdot\psi\equiv0$ in $\Omega\times(0,\infty)$, and if finally there exists some nonnegative $\Phi\in\CSp{2}{[0,\infty)}$ with $\Phi'>0$ on $(0,\infty)$ such that $n$ is a global weak $\Phi$--supersolution of \eqref{eq:system_very_weak_supersol} in the sense of Definition \ref{def:weak_super_sol}.
\end{definition}

If, on the other hand, $m+2\alpha>\frac{5}{3}$ we will obtain global weak solutions in the standard sense. Let us formulate this well-established concept for the sake of completeness in the following definition.

\begin{definition}\label{def:weak-solution}
Let $S\in\CSp{2}{\bomega\times[0,\infty)^2;\R^{3\times3}}$ satisfy \eqref{Sprop} for some $S_0>0$ and $\alpha\geq0$. A triple $(n,c,u)$ of functions
\begin{align*}
n&\in\LSploc{1}{\bomega\times[0,\infty)},\\
c&\in\LSploc{1}{[0,\infty);\W[1,2]},\\
u&\in L_{loc}^1\big([0,\infty); W_{0}^{1,1}\!\left(\Omega;\R^3\right)\!\big),
\end{align*}
satisfying $n\geq0$ and $c\geq0$ in $\bomega\times[0,\infty)$, $n\in \LSploc{1}{[0,\infty);\W[1,1]}$ and $cu\in\LSploc{1}{\bomega\times[0,\infty);\R^3}$, as well as $u\otimes u\in\LSploc{1}{\bomega\times[0,\infty);\R^{3\times3}}$ will be called a global weak solution of \eqref{CTndS}--\,\eqref{IC}, if $\nabla\cdot u=0$ in $\mathcal{D}'\big((\Omega\times(0,\infty)\big)$, if
\begin{align*}
n^{m-1}\nabla n\quad\text{and}\quad n\nabla c,\quad\text{as well as}\quad nu&\quad\text{belong to }\LSploc{1}{\bomega\times[0,\infty);\R^3},
\end{align*}
if equality \eqref{eq:very_weak_sol1} holds for all $\varphi\in\LSp{\infty}{\Omega\times(0,\infty)}\cap\LSp{2}{(0,\infty);\W[1,2]}$ with $\varphi_t\in\LSp{2}{\Omega\times(0,\infty)}$, which are compactly supported in $\bomega\times[0,\infty)$, if  \eqref{eq:very_weak_sol2} is fulfilled for all $\psi\in C_0^\infty\left(\Omega\times[0,\infty);\R^3\right)$ with $\nabla\cdot\psi\equiv0$ in $\Omega\times(0,\infty)$, and if finally for each $\varphi\in C_0^\infty\big(\bomega\times[0,\infty)\big)$ with $\frac{\partial\varphi}{\partial\nu}=0$ on $\romega\times(0,\infty)$, the equality
\begin{align}\label{eq:weak-sol-n}
-\intinfomega n\varphi_t&-\intomega n_0\,\varphi(\cdot,0)\nonumber\\&=-m\intinfomega n^{m-1} \big(\nabla n\cdot\nabla\varphi\big)+\intinfomega n\big(S(x,n,c)\nabla c\cdot\nabla\varphi\big)+\intinfomega n(u\cdot\nabla\varphi)
\end{align}
is satisfied.
\end{definition}

\begin{remark} i) If \eqref{eq:very_weak_supersol} is satisfied for $\Phi(s)\equiv s$ with equality, then $(n,c,u)$ is a global weak solution of \eqref{CTndS} in the sense of Definition \ref{def:weak-solution}, which shows that every global weak solution is also a global very weak solution.

ii) If the global very weak solution $(n,c,u)$ satisfies the regularity properties $n,c\in\CSp{0}{\bomega\times[0,\infty)}\cap\CSp{2,1}{\bomega\times(0,\infty)}$ and $u\in\CSp{0}{\bomega\times[0,\infty);\R^3}\cap\CSp{2,1}{\bomega\times(0,\infty);\R^3}$, it can be checked that the solution is also a global classical solution, i.e. one can find $P\in\CSp{1,0}{\bomega\times(0,\infty)}$ such that $(n,c,u,P)$ solves \eqref{CTndS} in the classical sense. See \cite[Lemma 2.1]{win15_chemorot} for the arguments involved. 
\end{remark}

\setcounter{equation}{0} 
\section{A family of regularized problems}\label{sec3:approxprob}
As a first step in the construction of global solutions in either of the senses above we will first adapt the approaches undertaken in \cite{win15_chemorot,Wang17-globweak-ksns,TB2017_nonlindiff} to our setting in order to approximate the system \eqref{CTndS} by problems in which the no-flux boundary condition of the first component reduces to a homogeneous Neumann boundary condition and which are solvable globally in time. With a family $(\rho_\epsi)_{\epsi\in(0,1)}\subset C_0^\infty(\Omega)$ of cut-off functions in $\Omega$ satisfying 
\begin{align*}
0\leq \rho_\epsi(x)\leq 1\text{ for all }x\in\Omega\quad\text{such that}\quad\rho_\epsi\nearrow1\text{ as }\epsi\searrow0,
\end{align*}
we define 
\begin{align}\label{Sepsi}
S_\epsi(x,n,c):=\rho_\epsi(x)S(x,n,c),\quad (x,n,c)\in\Omega\times[0,\infty)^2
\end{align}
and accordingly for $\epsi\in(0,1)$ consider regularized problems of the form
\begin{align}\label{approxprob}
\left\{
\begin{array}{r@{\,}c@{\, }c@{\,}l@{\quad}l@{\quad}l@{\,}c}
n_{\epsi t}&+&u_\epsi\cdot\!\nabla n_\epsi&=\nabla\cdot\big(m(n_\epsi+\epsi)^{m-1}\nabla n_\epsi-\frac{n_\epsi S_\epsi(x,n_\epsi,c_\epsi)}{(1+\epsi n_\epsi)^3}\nabla c_\epsi\big),\ &x\in\Omega,& t>0,\\
c_{\epsi t}&+&u_\epsi \cdot\!\nabla c_\epsi&=\Delta c_\epsi-c_\epsi+n_\epsi,\ &x\in\Omega,& t>0,\\
u_{\epsi t}&+&(Y_\epsi u_\epsi\cdot\nabla)u_\epsi&=\Delta u_\epsi+\nabla P_\epsi+n_\epsi\nabla\phi,\ &x\in\Omega,& t>0,\\
&&\quad\nabla\cdot u_\epsi&=0,\ &x\in\Omega,& t>0,\\
&&\quad\ \partial_\nu n_\epsi=\partial_\nu c_\epsi&=0,\qquad\qquad\quad u_\epsi=0,\ &x\in\romega,& t>0,\\
&&\qquad n_\epsi(x,0)&=n_0(x),\quad c_\epsi(x,0)=c_0(x),\quad u_\epsi(x,0)=u_0(x),\ &x\in\Omega,&
\end{array}\right.
\end{align}
where the Yosida approximation of the Stokes operator $Y_\epsi$ is given by
\begin{align*}
Y_\epsi\varphi:=(1+\epsi A)^{-1}\varphi\quad\text{for }\epsi\in(0,1)\text{ and }\varphi\in L_\sigma^2(\Omega).
\end{align*}
\subsection{Global existence of approximating solutions and basic properties}\label{sec31:locex}
By standard arguments involving well-established testing procedures and a Moser-type iteration one can readily verify that for all $m\geq1$ and $\alpha\geq0$ the classical solutions to the approximating system above are in fact global solutions, which in addition satisfy certain $\Lo[1]$--estimates.

\begin{lemma}\label{lem:glob_ex}
Let $\Omega\subset\R^3$ be a bounded domain with smooth boundary, $\phi\in\W[2,\infty]$, $\vartheta>3$, $m\geq1$ and $\alpha\geq0$. Suppose that $S\in\CSp{2}{\bomega\times[0,\infty)^2;\R^{3\times3}}$ satisfies \eqref{Sprop} for some $S_0>0$ and assume that $n_0,c_0$ and $u_0$ comply with \eqref{IR}. Then for any $\epsi\in(0,1)$, there exists a uniquely determined triple $(n_\epsi,c_\epsi,u_\epsi)$ of functions satisfying
\begin{align*}
n_\epsi&\in \CSp{0}{\bomega\times[0,\infty)}\cap \CSp{2,1}{\bomega\times(0,\infty)}, \\
c_\epsi&\in \CSp{0}{\bomega\times[0,\infty)}\cap \CSp{2,1}{\bomega\times(0,\infty)}\cap\CSp{0}{[0,\infty);\W[1,\vartheta]},\\
u_\epsi&\in \CSp{0}{\bomega\times[0,\infty);\R^3}\cap \CSp{2,1}{\bomega\times(0,\infty);\R^3},
\end{align*}
which, together with some $P_\epsi\in \CSp{1,0}{\bomega\times(0,\infty)}$, solve \eqref{approxprob} in the classical sense and fulfill $n_\epsi\geq0$ and $c_\epsi\geq0$ in $\bomega\times[0,\Tme)$, as well as
\begin{align}\label{eq:mass-cons-n}
\intomega n_\epsi(\cdot,t)=\intomega n_0\quad\text{for all }t\in(0,\infty)
\end{align}
and
\begin{align}\label{eq:l1-bound-c}
\intomega c_\epsi(\cdot,t)\leq \max\left\{\intomega n_0,\intomega c_0\right\}\quad\text{for all }t\in(0,\infty).
\end{align}
\end{lemma}

\begin{bew}
In light of the fact that $|S_\epsi(x,n,c)|\leq S_0$ on $\Omega\times[0,\infty)^2$ this is essentially already contained in \cite[Lemmata 3.1--\,3.4]{TB2017_nonlindiff} with minimal necessary adjustments. Let us briefly state the main ideas. The proof of local-in-time classical solutions on $\bomega\times(0,\Tme)$, where $\Tme\in(0,\infty]$ denotes the maximal existence time, can be achieved by adapting standard fixed point arguments as illustrated for similar chemotaxis frameworks in e.g. \cite[Lemma 2.1]{tao_winkler_chemohapto11-siam11}, \cite[Lemma 2.2]{Lan17-LocBddGlobSolNonlinDiff-JDE} and \cite[Lemma 2.1]{win_fluid_final}. The nonnegativity of $n_\epsi$ and $c_\epsi$ can then be established by relying on the maximum principle, whereas the $L^1$-regularity of $n_\epsi$ and $c_\epsi$ follows immediately from integrating the corresponding equations and, for $c_\epsi$, employment of an ODE comparison argument. To verify that the solution is indeed global in time we first rely on standard testing procedures to obtain that for fixed $\epsi\in(0,1)$ and $T\in(0,\Tme]$ with $T<\infty$ there exists $C_1>0$ such that
\begin{align*}
\intomega n_\epsi^6(\cdot,t)+\intomega c_\epsi^6(\cdot,t)\leq C_1\quad\text{holds for all }t\in(0,T).
\end{align*}
Relying on further testing procedures for the third equation (see also \cite[Lemma 3.9]{win_globweak3d-AHPN16}) and the smoothing properties of the Stokes operator (e.g. \cite[Lemma 3.1]{Win-ct_fluid_3d-CPDE15}) we find that for $\beta\in(\tfrac34,1)$ there exists $C_2(T)>0$ such that $\|A^\beta u_\epsi\|_{\Lo[2]}\leq C_2(T)$ for all $t\in(0,T)$. These bounds at hand we can go to testing the second equation by $-\Delta c_\epsi$ to first obtain $L^2$-information on $\nabla c_\epsi$, which, by standard semigroup estimates, can then be refined to a bound on $\|\nabla c_\epsi(\cdot,t)\|_{\Lo[\frac{11}2]}$ for all $t\in(0,T)$. Combining this with our previous bounds we can employ a Moser-type iteration (see e.g. \cite[Lemma A.1]{TaoWin-quasilinear_JDE12}) to finally conclude that in fact $\Tme=\infty$.
\end{bew}

\setcounter{equation}{0} 
\section{A priori estimates}\label{sec4:reg-est}
As our main focus will be on values $m\geq1$ and $\alpha\geq0$ which are both as small as possible, our main task will be to obtain regularity information independent on $\epsi\in(0,1)$, which restrict $m$ and $\alpha$ in the least possible way. As in particular no energy-structure is present in \eqref{approxprob} we are thereby task with finding a testing procedure, which captures as optimal conditions on these parameters as possible. Even obtaining an $L^2$--estimate for $n_\epsi$ seems to be far out of reach without gravely restricting either $m$ or $\alpha$. Thus, similar to the approach in \cite{TB2017_nonlindiff}, we decide to investigate a functional which for small values of $m$ and $\alpha$ is of sublinear growth, hoping to obtain a spatio-temporal bound on the gradient of $n_\epsi$, which we can refine later on to a regularity estimate beyond the $L^1$--estimate of Lemma \ref{lem:glob_ex}.
\subsection{Estimates capturing optimal conditions on \texorpdfstring{$m$ and $\alpha$}{m and alpha}}\label{sec41:regesti}
Let us start with an elementary identity laying the groundwork to impending testing procedures.
\begin{lemma}\label{lem:testing-nepsi-m+2a-1}
Let $m\geq1$, $\alpha\geq0$, $\beta\geq1$ be such that $m+\frac{\beta}{2}\alpha>1$, assume that $n_0,c_0$ and $u_0$ comply with \eqref{IR} and that $S\in\CSp{2}{\bomega\times[0,\infty)^2;\R^{3\times3}}$. Then for any $\epsi\in(0,1)$ and each $\varphi\in C^\infty\left(\bomega\times[0,\infty)\right)$ with $\frac{\partial\varphi}{\partial\nu}=0$ on $\romega\times(0,\infty)$ the classical solution $(n_\epsi,c_\epsi,u_\epsi)$ of \eqref{approxprob} satisfies
\begin{align}\label{eq:testing-nepsi}
\frac{\intd}{\intd t}\intomega& (n_\epsi+\epsi)^{m+\beta\alpha-1}\nonumber\\=\ &\frac{m(m+\beta\alpha-1)(2-(m+\beta\alpha))}{(m+\frac\beta2\alpha-1)^2}\intomega \big|\nabla(n_\epsi+\epsi)^{m+\frac\beta2\alpha-1}\big|^2\\
&-\frac{(m+\beta\alpha-1)(2-(m+\beta\alpha))}{m+\frac\beta2\alpha-1}\intomega\frac{n_\epsi(n_\epsi+\epsi)^{\frac\beta2\alpha-1}}{(1+\epsi n_\epsi)^3}\big(\nabla (n_\epsi+\epsi)^{m+\frac\beta2\alpha-1}\cdot S_\epsi(x,n_\epsi,c_\epsi)\nabla c_\epsi\big)\nonumber
%
\end{align}
 on $(0,\infty)$.
\end{lemma} 

\begin{bew}
Drawing on the first equation of \eqref{approxprob} straightforward calculations show that
\begin{align*}
&\frac{\intd}{\intd t}\intomega (n_\epsi+\epsi)^{m+\beta\alpha-1}\\
=\ &(m+\beta\alpha-1)\intomega(n_\epsi+\epsi)^{m+\beta\alpha-2}\nabla\cdot\Big(m(n_\epsi+\epsi)^{m-1}\nabla n_\epsi-\frac{n_\epsi}{(1+\epsi n_\epsi)^3}S_\epsi(\cdot,n_\epsi,c_\epsi)\nabla c_\epsi\Big)\\&\qquad-(m+\beta\alpha-1)\intomega(n_\epsi+\epsi)^{m+\beta\alpha-2}(u_\epsi\cdot\nabla n_\epsi)
\end{align*}
holds on $(0,\infty)$. Making use of the fact that $\nabla\cdot u_\epsi\equiv 0$ in $\Omega\times(0,\infty)$ as well as the imposed boundary conditions, we find that upon integration by parts and appropriate reformulation of some terms the asserted equality follows immediately.
\end{bew}

Depending on the sign of $2-(m+2\alpha)$, we will multiply the equality of Lemma \ref{lem:testing-nepsi-m+2a-1} with either positive or negative constants and then estimate. Combining the resulting inequality with a standard testing procedure for the second equation we will derive some information on $(n_\epsi+\epsi)^{m+2\alpha-1}$, $\nabla(n+\epsi)^{m+\alpha-1}$, $c_\epsi^2$ and $\nabla c_\epsi^2$. This approach has been undertaken previously in e.g. \cite[Lemma 4.1]{Wang17-globweak-ksns} and \cite[Lemma 4.2]{TB2017_nonlindiff}.

\begin{lemma}\label{lem:bounds}
Let $m\geq1$, $\alpha\geq0$ be such that $m+\alpha>\frac{4}{3}$, suppose that $n_0,c_0$ and $u_0$ fulfill \eqref{IR} and assume that $S\in\CSp{2}{\bomega\times[0,\infty)^2;\R^{3\times 3}}$ satisfies \eqref{Sprop} with some $S_0>0$. Then there exists some $C>0$ such that for all $\epsi\in(0,1)$ the global classical solution $(n_\epsi,c_\epsi,u_\epsi)$ of \eqref{approxprob} satisfies
\begin{align}\label{eq:lem-bounds-eq}
\intomega (n_\epsi+\epsi)^{m+2\alpha-1}(\cdot,t)+\intomega c_\epsi^2(\cdot,t) +\int_t^{t+1}\!\intomega\big|\nabla (n_\epsi+\epsi)^{m+\alpha-1}\big|^2+\int_t^{t+1}\!\intomega|\nabla c_\epsi|^2\leq C
\end{align}
for all $t\geq0$.
\end{lemma}

\begin{bew}
Since the main part of the procedure does not differ to greatly from the setting with a scalar sensitivity as discussed in \cite[Lemma 4.2]{TB2017_nonlindiff}, we will only cover the main ideas. First assume $m+2\alpha<2$. Employing Lemma \ref{lem:testing-nepsi-m+2a-1} with $\beta=2$ and multiplying the equality by $-\frac{1}{(m+2\alpha-1)}$ we can make use of Young's inequality and the fact that $\Big|\frac{n_\epsi(n_\epsi+\epsi)^{\alpha-1}S_\epsi(\cdot,n_\epsi,c_\epsi)}{(1+\epsi n_\epsi)^3}\Big|\leq S_0$ in $\Omega\times(0,\infty)$ to find that
\begin{align}\label{eq:bounds-test-n}
-\frac{1}{m+2\alpha-1}&\frac{\intd}{\intd t}\intomega(n_\epsi+\epsi)^{m+2\alpha-1}\\&\leq-\frac{m(2-(m+2\alpha))}{2(m+\alpha-1)^2}\!\intomega\big|\nabla(n_\epsi+\epsi)^{m+\alpha-1}\big|^2+\frac{S_0^2(2-(m+2\alpha))}{2m}\!\intomega|\nabla c_\epsi|^2\nonumber
\end{align}
in $(0,\infty)$ for all $\epsi\in(0,1)$. Testing the second equation of \eqref{approxprob} by $c_\epsi$, we find that an application of Hölder's inequality and the embedding $\W[1,2]\hookrightarrow\Lo[6]$ entail the existence of $C_1>0$ such that
\begin{align}\label{eq:bounds-test-c}
\frac{\intd}{\intd t}\intomega c_\epsi^2(\cdot,t)+\intomega|\nabla c_\epsi(\cdot,t)|^2+\intomega c_\epsi^2(\cdot,t)\leq C_1\|n_\epsi(\cdot,t)\|_{\Lo[\nfrac{6}{5}]}^2
\end{align}
for all $t>0$ and all $\epsi\in(0,1)$, where we used that $u_\epsi$ is a solenoidal vector field. Moreover, drawing on the \GNI, the nonnegativity of $n_\epsi$, the mass conservation featured in Lemma \ref{lem:glob_ex} and the fact that $\epsi\in(0,1)$, we find $C_2>0$ such that
\begin{align*}
C_1\|n_\epsi\|_{\Lo[\nfrac{6}{5}]}^2\leq C_1\big\|(n_\epsi+\epsi)^{m+\alpha-1}\big\|_{\Lo[\frac{6}{5(m+\alpha-1)}]}^\frac{2}{m+\alpha-1}\leq C_2\big\|\nabla(n_\epsi+\epsi)^{m+\alpha-1}\big\|_{\Lo[2]}^{\frac{2}{6(m+\alpha)-7}}+C_2
\end{align*}
holds on $(0,\infty)$ for all $\epsi\in(0,1)$, and, since $m+\alpha>\frac{4}{3}$ implies $\frac{2}{6(m+\alpha)-7}<2$, an application of Young's inequality thereby provides $C_3>0$ such that
\begin{align*}
C_1\|n_\epsi(\cdot,t)\|_{\Lo[\nfrac{6}{5}]}^2\leq \frac{m^2}{4S_0^2(m+\alpha-1)^2}\intomega\big|\nabla(n_\epsi+\epsi)^{m+\alpha-1}(\cdot,t)\big|^2+C_3\quad\text{for all }t>0\text{ and all }\epsi\in(0,1).
\end{align*}
Thus, combining \eqref{eq:bounds-test-n} with a multiple of \eqref{eq:bounds-test-c} and the estimate above we have
\begin{align}\label{eq:bounds-ode}
y_\epsi'(t)+y_\epsi(t)+g_\epsi(t)\leq C_4\quad\text{for all }t>0\text{ and all }\epsi\in(0,1),
\end{align}
where we have set $C_4:=\frac{C_3 S_0^2(2-(m+2\alpha))}{m}>0$,
\begin{align*}
y_\epsi(t):=-\frac{1}{m+2\alpha-1}\intomega(n_\epsi+\epsi)^{m+2\alpha-1}(\cdot,t)+\frac{S_0^2(2-(m+2\alpha))}{m}\intomega c_\epsi^2(\cdot,t),\quad t>0,
\end{align*}
and
\begin{align*}
g_\epsi(t):=\frac{m(2-(m+2\alpha))}{4(m+\alpha-1)^2}\intomega\big|(n_\epsi+\epsi)^{m+\alpha-1}(\cdot,t)\big|^2+\frac{S_0^2(2-(m+2\alpha))}{2m}\intomega|\nabla c_\epsi(\cdot,t)|^2,\quad t>0.
\end{align*}
An ODE comparison implies the existence of $C_5>0$ satisfying $y_\epsi(t)\leq C_5$ for all $t>0$ and all $\epsi\in(0,1)$, which together with the definition of $y_\epsi$, the positivity of $\intomega c_\epsi^2$, the fact that $m+2\alpha-1<1$ and  Lemma \ref{lem:glob_ex} shows that for some $C_6>0$ we have $|y_\epsi(t)|\leq C_6$ for all $t>0$ and all $\epsi\in(0,1)$, proving bounds for the first two summands in \eqref{eq:lem-bounds-eq}. For the remaining terms we integrate \eqref{eq:bounds-ode} with respect to time to find that
\begin{align*}
\int_t^{t+1}\!g_\epsi(s)\intd s\leq y_\epsi(t)-y_\epsi(t+1)-\int_t^{t+1}\!y_\epsi(s)\intd s+C_4
\end{align*}
holds for all $t>0$ and all $\epsi\in(0,1)$. Hence, the previously discussed boundedness of $|y_\epsi(t)|$ entails the boundedness of the latter two terms in \eqref{eq:lem-bounds-eq}.

In the case of $m+2\alpha>2$ we can follow the same arguments as above with multiplying the equation of Lemma \ref{lem:testing-nepsi-m+2a-1} this time with $\frac{1}{m+2\alpha-1}$ to obtain a similar ODE to \eqref{eq:bounds-ode}, which then allows us to conclude the asserted bounds in similar fashion. If $m+2\alpha=2$, we note that $m+\alpha-1=1-\alpha$ and that moreover $\alpha\leq\frac{1}{2}$ due to $m\geq1$. Thus, estimating
\begin{align*}
\frac{\intd}{\intd t}\intomega(n_\epsi\ln n_\epsi)(\cdot,t)\leq-\frac{m}{2(1-\alpha)^2}\intomega\big|\nabla(n_\epsi+\epsi)^{1-\alpha}(\cdot,t)\big|^2+\frac{S_0^2}{2m}\intomega|\nabla c_\epsi(\cdot,t)|^2,
\end{align*} 
for all $t>0$ and all $\epsi\in(0,1)$, and combining with \eqref{eq:bounds-test-c} we obtain an inequality of the form
\begin{align*}
\frac{\intd}{\intd t}\Big(\intomega(n_\epsi\ln n_\epsi)(\cdot,t)+\intomega c_\epsi^2(\cdot,t)\Big)+C_7\intomega\big|\nabla (n_\epsi+\epsi)^{1-\alpha}(\cdot,t)\big|^2+C_7\intomega|\nabla c_\epsi(\cdot,t)|^2+C_7\intomega c_\epsi^2(\cdot,t)\leq C_8,
\end{align*}
with some $C_7>0$ and $C_8>0$. By means of the \GNI\ and the evident estimate $x\ln x\leq x^\frac{5}{3}$ for $x>0$ we have $C_9>0$ satisfying
\begin{align*}
\intomega n_\epsi\ln n_\epsi\leq\|(n_\epsi+\epsi)^{1-\alpha}\|_{\Lo[\frac{5}{3(1-\alpha)}]}^\frac{5}{3(1-\alpha)}\leq C_9\|\nabla (n_\epsi+\epsi)^{(1-\alpha)}\|_{\Lo[2]}^\frac{4}{5-6\alpha}+C_9\quad\text{on}\quad(0,\infty).
\end{align*}
Because of $\frac{4}{5-6\alpha}\leq2$ for $\alpha\leq\frac12$, this now allows to pursue a similar reasoning as before, while making use of the fact that $s\ln s\geq-\frac{1}{e}$ for all $s>0$.
\end{bew} 

While the main idea of utilizing the latter spatio-temporal bound for $\nabla (n_\epsi+\epsi)^{m+\alpha-1}$ to establish time-space bounds for $n_\epsi+\epsi$ remains unchanged from the previous works \cite[Lemma 4.2]{Wang17-globweak-ksns} and \cite[Lemma 4.3]{TB2017_nonlindiff}, we have to treat the term more delicate in order to prepare sufficient information for the limiting procedure later on.

\begin{lemma}\label{lem:st-bound-nepsi}
Let $m\geq1$, $\alpha\geq0$ be such that $m+\alpha>\frac{4}{3}$ and assume that $n_0,c_0$ and $u_0$ comply with \eqref{IR} and that $S\in\CSp{2}{\bomega\times[0,\infty)^2;\R^{3\times 3}}$ fulfills \eqref{Sprop} with some $S_0>0$. Then for all $p\in\big(1,6(m+\alpha-1)\big)$ there exists $C>0$ such that for all $\epsi\in(0,1)$ the solution $(n_\epsi,c_\epsi,u_\epsi)$ of \eqref{approxprob} satisfies
\begin{align}\label{eq:st-bound-np}
\int_t^{t+1}\! \big\|n_\epsi(\cdot,s)+\epsi\big\|_{\Lo[p]}^{\frac{2p(m+\alpha-\frac{7}{6})}{p-1}}\intd s\leq C\quad\text{for all }t\geq 0.
\end{align}
In particular, there exist $r\in(1,2)$ and $C>0$ such that
\begin{align}\label{eq:st-bound-specialcases}
\int_t^{t+1}\!\big\|n_\epsi(\cdot,s)+\epsi\big\|_{\Lo[\frac{6r}{6-r}]}^{\frac{2r}{2-r}}\intd s\leq C\quad\text{and}\quad\int_t^{t+1}\big\|n_\epsi(\cdot,s)+\epsi\big\|_{\Lo[2(m+\alpha)-\frac{4}{3}]}^{2(m+\alpha)-\frac{4}{3}}\intd s\leq C
\end{align}
hold for each $\epsi\in(0,1)$ and all $t\geq0$.
\end{lemma}

\begin{bew} We employ reasoning similar to \cite[Lemma 4.2]{Wang17-globweak-ksns} and \cite[Lemma 4.3]{TB2017_nonlindiff}.
Due to $p\in\big(1,6(m+\alpha-1)\big)$ and $m+\alpha>\frac43>\frac76$ we can utilize the \GNI\ (e.g. \cite[Lemma 2.3]{lankchapto15}) to find $C_1>0$ such that with
\begin{align*}
a=\frac{m+\alpha-1-\frac{m+\alpha-1}{p}}{m+\alpha-1+\frac{1}{3}-\frac{1}{2}}=\frac{p-1}{p}\cdot\frac{6(m+\alpha-1)}{6m+6\alpha-7}\in(0,1)
\end{align*}
the inequality
\begin{align*}
\int_t^{t+1}\!\!\!\big\|n_\epsi&(\cdot,s)+\epsi\big\|_{\Lo[p]}^{\frac{2p(m+\alpha-\frac{7}{6})}{p-1}}\intd s=\int_t^{t+1}\!\!\!\big\|(n_\epsi+\epsi)^{m+\alpha-1}(\cdot,s)\big\|_{\Lo[\frac{p}{m+\alpha-1}]}^{\frac{2p}{p-1}\cdot\frac{6m+6\alpha-7}{6(m+\alpha-1)}}\intd s\\
&\leq C_1\!\int_t^{t+1}\!\!\!\big\|\nabla(n_\epsi+\epsi)^{m+\alpha-1}(\cdot,s)\big\|_{\Lo[2]}^{\frac{2p}{p-1}\cdot\frac{6m+6\alpha-7}{6(m+\alpha-1)}\cdot a}\big\|(n_\epsi+\epsi)^{m+\alpha-1}(\cdot,s)\big\|_{\Lo[\frac{1}{m+\alpha-1}]}^{\frac{2p}{p-1}\cdot\frac{6m+6\alpha-7}{6(m+\alpha-1)}\cdot(1-a)}\intd s\\&\hspace*{6.3cm}+C_1\!\int_t^{t+1}\!\!\!\big\|(n_\epsi+\epsi)^{m+\alpha-1}(\cdot,s)\big\|_{\Lo[\frac{1}{m+\alpha-1}]}^{\frac{2p}{p-1}\cdot\frac{6m+6\alpha-7}{6(m+\alpha-1)}}\intd s
\end{align*}
holds for all $t\geq0$ and all $\epsi\in(0,1)$. Combined with the mass conservation of $n_\epsi$, as established in Lemma \ref{lem:glob_ex}, this implies the existence of $C_2>0$ such that
\begin{align*}
\int_t^{t+1}\!\!\!\big\|n_\epsi&(\cdot,s)+\epsi\big\|_{\Lo[p]}^{\frac{2p(m+\alpha-\frac{7}{6})}{p-1}}\intd s\leq C_2
\int_t^{t+1}\!\!\!\big\|\nabla(n_\epsi+\epsi)^{m+\alpha-1}(\cdot,s)\big\|_{\Lo[2]}^{2}	\intd s+C_2
\end{align*}
holds for all $t\geq0$ and all $\epsi\in(0,1)$, which proves \eqref{eq:st-bound-np} under consideration of Lemma \ref{lem:bounds}. For the first special case in \eqref{eq:st-bound-specialcases} we first note that due to $m+\alpha>\frac{4}{3}$ the interval $I:=\big(1,\min\{\frac{6(m+\alpha-1)}{m+\alpha},2\}\big)$ is not empty and that for $r\in I$ we have $q:=\frac{6r}{6-r}\in\big(1,6(m+\alpha-1)\big)$.
Moreover, $r<\frac{6(m+\alpha-1)}{m+\alpha}$ together with $m+\alpha>\frac{7}{6}$ also readily implies $r<\frac{1+2(m+\alpha-\frac76)}{m+\alpha}$ and hence
\begin{align*}
\frac{2q(m+\alpha-\frac76)}{q-1}=\frac{12r(m+\alpha-\frac76)}{7r-6}>\frac{2r}{2-r}.
\end{align*}
Thus, the first special case follows from \eqref{eq:st-bound-np} with $p=\frac{6r}{6-r}$. For the second bound in \eqref{eq:st-bound-specialcases} we work along similar lines noting that, again due to $m+\alpha>\frac43$, $2(m+\alpha)-\frac{4}{3}\in\big(1,6(m+\alpha-1)\big)$ and that $2(m+\alpha)-\frac{4}{3}=\frac{2(m+\alpha-\frac76)(2(m+\alpha)-\frac43)}{2(m+\alpha)-\frac43-1}$, making the first part of the lemma applicable once more.
\end{bew}

Let us also briefly establish some supplementary spatio-temporal estimates under the additional assumption that $m+\alpha\leq 2$. These bounds follow in a straightforward fashion from Lemma \ref{lem:bounds} and Lemma \ref{lem:st-bound-nepsi}, and will later form a cornerstone in obtaining convergence properties necessary to pass to the limit in the integrals making up the global weak $\Phi$-supersolution for $\Phi(s)=(s+1)^{m+2\alpha-1}$.

\begin{corollary}\label{cor:1}
Let $m\geq1$, $\alpha\geq0$ be such that $\tfrac43<m+\alpha\leq 2$, suppose that $n_0,c_0$ and $u_0$ fulfill \eqref{IR} and assume that $S\in\CSp{2}{\bomega\times[0,\infty)^2;\R^{3\times 3}}$ satisfies \eqref{Sprop} with some $S_0>0$. Then there exists some $C_1>0$ such that for all $\epsi\in(0,1)$ the global classical solution $(n_\epsi,c_\epsi,u_\epsi)$ of \eqref{approxprob} satisfies
\begin{align}\label{eq:cor1}
\int_t^{t+1}\!\intomega\big|\nabla(n_\epsi+1)^{m+\alpha-1}\big|^2+\int_t^{t+1}\!\intomega\big|(n_\epsi+1)^{\frac{m+2\alpha-3}{2}}(n_\epsi+\epsi)^{\frac{m-1}{2}}\nabla n_\epsi\big|^2\leq C_1,
\end{align}
for all $t\geq0$. Moreover, there exist $p>2$, $r>1$ and $C_2>0$ such that
\begin{align}\label{eq:cor2}
\int_t^{t+1}\big\|(n_\epsi+1)^{m+\alpha-1}\big\|_{\Lo[p]}^p&\leq C_2\quad\text{and}\quad\int_t^{t+1}\big\|(n_\epsi+1)^\alpha(n_\epsi+\epsi)^{m-1}\big\|_{\Lo[p]}^p\leq C_2,
\intertext{as well as}\label{eq:cor3}
&\int_t^{t+1}\big\|(n_\epsi+1)^{m+2\alpha-1}\big\|_{\Lo[\frac{6r}{6-r}]}^{\frac{2r}{2-r}}\leq C_2
\end{align}
hold for each $\epsi\in(0,1)$ and all $t\geq0$.
\end{corollary}

\begin{bew}
Due to $m+\alpha\in(\frac43,2]$ it is obvious that
\begin{align*}
\int_t^{t+1}\!\intomega\big|\nabla(n_\epsi+1)^{m+\alpha-1}\big|^2&=\int_t^{t+1}\!\intomega(n_\epsi+1)^{2(m+\alpha-2)}|\nabla n_\epsi|^2\\&\leq\int_t^{t+1}\!\intomega(n_\epsi+\epsi)^{2(m+\alpha-2)}|\nabla n_\epsi|^2=\int_t^{t+1}\!\intomega\big|\nabla(n_\epsi+\epsi)^{m+\alpha-1}\big|^2
\end{align*}
holds for all $\epsi\in(0,1)$ and $t\geq0$, whereupon the boundedness of the first term in \eqref{eq:cor1} follows immediately from Lemma \ref{lem:bounds}. The bound for the second term contained in \eqref{eq:cor1} then is a direct consequence of the first bound in light of the fact that $m\geq1$.
Reiterating the proof of Lemma \ref{lem:st-bound-nepsi} for $(n_\epsi+1)$ instead of $(n_\epsi+\epsi)$, while relying on \eqref{eq:cor1}, we find that for all $q\in\big(1,6(m+\alpha-1)\big)$ there exists $C>0$ such that
\begin{align*}
\int_t^{t+1}\big\|n_\epsi(\cdot,s)+1\big\|_{\Lo[q]}^{\frac{2q(m+\alpha-\frac76)}{q-1}}\intd s\leq C\quad\text{for all }t\geq0.
\end{align*}
This spatio-temporal estimate at hand, straightforward calculations, similar to those undertaken to prove the special cases presented in Lemma \ref{lem:st-bound-nepsi}, verify \eqref{eq:cor2} and \eqref{eq:cor3}, due to the facts that $m\geq1$, $\alpha\geq0$ and $m+\alpha>\frac43$.
\end{bew}
	
\subsection{Estimates involving the fluid component \texorpdfstring{$u_\epsi$}{}}

We will briefly state \cite[Lemma 3.4]{LL18-ClassSolLogCTSingSens} without proof. This result will be applied to a differential inequality for $\intomega|u_\epsi(\cdot,t)|^2$ in the lemma thereafter to obtain a first boundedness information on the fluid component, which can then be refined to additional spatio-temporal bounds.

\begin{lemma}\label{lem:diffineq-lemma}
For some $T\in(0,\infty]$ let $y\in\CSp{1}{(0,T)}\cap\CSp{0}{[0,T)}$, $h\in\CSp{0}{[0,T)}$, $C>0$, $a>0$ satisfy
\begin{align*}
y'(t)+ay(t)\leq h(t),\qquad\int_{(t-1)_+}^t h(s)\intd s\leq C
\end{align*}
for all $t\in(0,T)$. Then $y\leq y(0)+\frac{C}{1-e^{-a}}$ throughout $(0,T)$.
\end{lemma}

Drawing on Lemmata \ref{lem:st-bound-nepsi} and \ref{lem:diffineq-lemma}, as well as Hölder's inequality we are now in a position to utilize quite standard arguments, which and have been successfully employed before in e.g. \cite[Lemmata 3.5 and 3.6]{win_globweak3d-AHPN16} and \cite[Lemma 4.3]{Wang17-globweak-ksns}.

\begin{lemma}\label{lem:epsi-ind-est-u}
Let $m\geq1$, $\alpha\geq0$ be such that $m+\alpha>\frac{4}{3}$ and assume that $n_0,c_0$ and $u_0$ comply with \eqref{IR} and that $S\in\CSp{2}{\bomega\times[0,\infty)^2;\R^{3\times 3}}$ fulfills \eqref{Sprop} with some $S_0>0$. Then there exists $C>0$ such that for all $\epsi\in(0,1)$ the solution $(n_\epsi,c_\epsi,u_\epsi)$ of \eqref{approxprob} satisfies
\begin{align*}
\intomega |u_\epsi(\cdot,t)|^2+\int_t^{t+1}\!\intomega|\nabla u_\epsi|^2+\int_t^{t+1}\!\|u_\epsi\|_{\Lo[6]}^2\leq C
\end{align*}
for all $t\geq0$.
\end{lemma}

\begin{bew}
Multiplication of the third equation in \eqref{approxprob} by $u_\epsi$, integration by parts and an application of the Hölder inequality show that
\begin{align}\label{eq:epsi-ind-est-u-testing}
\frac{1}{2}\frac{\intd}{\intd t}\intomega|u_\epsi|^2(\cdot,t)+\intomega|\nabla u_\epsi(\cdot,t)|^2\leq\|\nabla\phi\|_{\Lo[\infty]}\|u_\epsi(\cdot,t)\|_{\Lo[6]}\|n_\epsi(\cdot,t)\|_{\Lo[\nfrac{6}{5}]}
\end{align}
holds for all $t>0$ and all $\epsi\in(0,1)$. Recalling the embedding $W_{0,\sigma}^{1,2}(\Omega)\hookrightarrow\Lo[6]$ and the Poincaré inequality we find $C_1>0$ satisfying
\begin{align}\label{eq:epsi-ind-est-u-pc+emb}
\|u_\epsi(\cdot,t)\|_{\Lo[6]}^2\leq C_1\intomega|\nabla u_\epsi(\cdot,t)|^2\quad\text{for all }t>0\text{ and all }\epsi\in(0,1),
\end{align}
which upon combination with \eqref{eq:epsi-ind-est-u-testing}, \eqref{phidef} and Young's inequality entails the existence of $C_2>0$ such that
\begin{align}\label{eq:epsi-ind-est-u-testing2}
\frac{1}{2}\frac{\intd}{\intd t}\intomega|u_\epsi|^2(\cdot,t)+\frac{1}{2}\intomega|\nabla u_\epsi(\cdot,t)|^2\leq C_2\|n_\epsi(\cdot,t)\|_{\Lo[\nfrac{6}{5}]}^2
\end{align}
is valid for all $t>0$ and all $\epsi\in(0,1)$. Due to Lemma \ref{lem:st-bound-nepsi} implying the existence of $C_3>0$ satisfying $\int_t^{t+1}\|n_\epsi(\cdot,t)\|_{\Lo[\nfrac{6}{5}]}^2\leq C_3$ for all $t>0$, we find that by estimating the gradient term by means of the Poincaré inequality from below and then employing Lemma \ref{lem:diffineq-lemma}, there exists $C_4>0$ such that
\begin{align*}
\intomega|u_\epsi|^2(\cdot,t)\leq C_4\quad\text{for all }t>0\text{ and all }\epsi\in(0,1).
\end{align*}
The estimate for $\intomega|u_\epsi|^2$ at hand, we can integrate \eqref{eq:epsi-ind-est-u-testing2} with respect to time to obtain that
\begin{align*}
\int_t^{t+1}\!\intomega|\nabla u_\epsi|^2\leq 2C_4+2C_2C_3\quad\text{for all }t>0\text{ and all }\epsi\in(0,1),
\end{align*}
which also immediately implies
\begin{align*}
\int_t^{t+1}\!\|u_\epsi\|_{\Lo[6]}^2\leq 2C_1C_4+2C_1C_2C_3\quad\text{for all }t>0\text{ and }\epsi\in(0,1),
\end{align*}
in light of \eqref{eq:epsi-ind-est-u-pc+emb}, and thus concludes the proof.
\end{bew}

With a first set of $\epsi$-independent estimates for the fluid component at hand, let us also briefly derive some spatio-temporal estimates for the combined quantities $n_\epsi u_\epsi$ and $(n_\epsi+1)^{m+2\alpha-1}u_\epsi$, which will be a cornerstone in treating the integrals appearing in the solution concepts which correspond to the convective term present in \eqref{approxprob}.

\begin{lemma}\label{lem:nepsi-uepsi-bounds}
Let $m\geq1$, $\alpha\geq0$ be such that $m+\alpha>\frac{4}{3}$ and assume that $n_0,c_0$ and $u_0$ comply with \eqref{IR} and that $S\in\CSp{2}{\bomega\times[0,\infty)^2;\R^{3\times 3}}$ fulfills \eqref{Sprop} with some $S_0>0$. Then there exist $r>1$ and $C_1>0$ such that for all $\epsi\in(0,1)$ the solution $(n_\epsi,c_\epsi,u_\epsi)$ of \eqref{approxprob} satisfies 
\begin{align*}
\int_t^{t+1}\!\intomega|n_\epsi u_\epsi|^r\leq C_1\quad\text{for all }t\geq0.
\end{align*}
If, moreover, $\frac43<m+\alpha\leq 2$, then there are $s>1$ and $C_2>0$ such that
\begin{align*}
\int_t^{t+1}\!\intomega\big|(n_\epsi+1)^{m+2\alpha-1}u_\epsi\big|^s\leq C_2
\end{align*}
hold for each $\epsi\in(0,1)$ and all $t\geq0$.
\end{lemma}

\begin{bew}
For any $r\in(1,2)$ an employment of the Hölder and Young inequalities to shows that
\begin{align*}
\int_t^{t+1}\!\intomega\big|n_\epsi u_\epsi\big|^r&\leq\int_t^{t+1}\big\|(n_\epsi+\epsi)u_\epsi\big\|_{\Lo[r]}^r\\
&\leq \int_t^{t+1}\|n_\epsi+\epsi\|_{\Lo[\frac{6r}{6-r}]}^r\|u_\epsi\|_{\Lo[6]}^r\leq\int_t^{t+1}\|n_\epsi+\epsi\|_{\Lo[\frac{6r}{6-r}]}^\frac{2r}{2-r}+\int_t^{t+1}\|u_\epsi\|_{\Lo[6]}^2
\end{align*}
holds for all $t\geq0$. Thus, taking $r>1$ as provided by Lemma \ref{lem:st-bound-nepsi}, the proof of the first assertion follows immediately from combining the estimate above with Lemmata \ref{lem:st-bound-nepsi} and \ref{lem:epsi-ind-est-u}. In a similar fashion we find that for $s\in(1,2)$ we have
\begin{align*}
\int_t^{t+1}\!\intomega\big|(n_\epsi+1)^{m+2\alpha-1}u_\epsi\big|^s\leq\int_t^{t+1}\big\|(n_\epsi+1)^{m+2\alpha-1}\big\|_{\Lo[\frac{6s}{6-s}]}^\frac{2s}{2-s}+\int_t^{t+1}\|u_\epsi\|_{\Lo[6]}^2
\end{align*}
for all $t\geq0$ and hence the second part of the Lemma is implied by Corollary \ref{cor:1} and Lemma \ref{lem:epsi-ind-est-u}. 
\end{bew}

\subsection{Time regularity}\label{sec42:timreg}

Having in mind an Aubin-Lions type argument to conclude the existence of limit objects of our approximate solution $(n_\epsi,c_\epsi,u_\epsi)$ when taking $\epsi\searrow0$, we still require regularity estimates for the time derivatives. Relying on the bounds established in the previous sections alone does not yet yield sufficient information on terms appearing in our estimation process.

\begin{lemma}\label{lem:nab-est-alphaover2}
Let $m\geq1$, $\alpha\geq0$ be such that $m+\alpha>\frac{4}{3}$, suppose that $n_0,c_0$ and $u_0$ comply with \eqref{IR} and assume that $S\in\CSp{2}{\bomega\times[0,\infty)^2;\R^{3\times3}}$ satisfies \eqref{Sprop} with some $S_0>0$. Then there exists $C>0$ such that for all $\epsi\in(0,1)$ the global classical solution $(n_\epsi,c_\epsi,u_\epsi)$ of \eqref{approxprob} satisfies
\begin{align*}
\int_t^{t+1}\!\intomega\big|\nabla(n_\epsi+\epsi)^{m+\frac{\alpha}{2}-1}\big|^2\leq C
\end{align*}
for all $t\geq0$.
\end{lemma}

\begin{bew}
Similar to the proof of Lemma \ref{lem:bounds} we first assume $m+\alpha<2$ and employ Lemma \ref{lem:testing-nepsi-m+2a-1} for $\beta=1$ and multiply the equality by $-\frac{1}{m+\alpha-1}$ to find that upon one application of Young's inequality that
\begin{align}\label{eq:nab-est-alphaover2-eq1}
-\frac{1}{m+\alpha-1}\frac{\intd}{\intd t}\intomega(n_\epsi+\epsi)^{m+\frac{\alpha}{2}-1}&+\frac{m(2-(m+\alpha))}{2(m+\frac{\alpha}{2}-1)^2}\intomega\big|\nabla(n_\epsi+\epsi)^{m+\frac{\alpha}{2}-1}\big|^2\\&\leq\frac{2-(m+\alpha)}{2m}\intomega\frac{n_\epsi^2(n_\epsi+\epsi)^{\alpha-2}}{(1+\epsi n_\epsi)^6}|S_\epsi(x,n_\epsi,c_\epsi)|^2|\nabla c_\epsi|^2\nonumber
\end{align}
holds on $(0,\infty)$. Noting that by $S_\epsi\leq S$ on $\Omega\times[0,\infty)^2$ and \eqref{Sprop} we have
\begin{align}\label{eq:nab-est-alphaover2-eq2}
\frac{n_\epsi^2(n_\epsi+\epsi)^{\alpha-2}}{(1+\epsi n_\epsi)^6}|S_\epsi(x,n_\epsi,c_\epsi)|^2\leq \frac{S_0^2(n_\epsi+\epsi)^{\alpha}}{(1+n_\epsi)^{2\alpha}}\leq S_0^2
\end{align}
we find upon integration of \eqref{eq:nab-est-alphaover2-eq1}, whilst also making use of the nonnegativity of $n_\epsi$ throughout $\bomega\times[0,\infty)$, that
\begin{align*}
\frac{m(2-(m+\alpha))}{2(m+\frac{\alpha}{2}-1)^2}&\int_t^{t+1}\!\intomega\big|\nabla(n_\epsi+\epsi)^{m+\frac{\alpha}{2}-1}\big|^2\\&\leq \frac{1}{m+\alpha-1}\intomega(n_\epsi+\epsi)^{m+\frac{\alpha}{2}-1}(\cdot,t+1)+\frac{S_0^2(2-(m+\alpha))}{2m}\int_t^{t+1}\!\intomega|\nabla c_\epsi|^2
\end{align*}
for all $t\geq0$ which proves the asserted bound for $m+\alpha<2$ in light of Lemma \ref{lem:bounds}. Identical arguments also work for $m+\alpha>2$ if one considers $\varphi=\frac{1}{m+\alpha-1}$. For $m+\alpha=2$ however, we will consider the time-evolution of $\intomega n_\epsi\ln n_\epsi$ to find that
\begin{align}\label{eq:nab-est-alphaover2-eq3}
\frac{\intd}{\intd t}\intomega n_\epsi\ln n_\epsi+\frac{m}{2(1-\frac{\alpha}{2})^2}\intomega\big|\nabla(n_\epsi+\epsi)^{1-\frac{\alpha}{2}}\big|^2\leq \frac{S_0^2}{2m}\intomega|\nabla c_\epsi|^2
\end{align}
on $(0,\infty)$.
where we used estimations akin to those in \eqref{eq:nab-est-alphaover2-eq2} and that $m-1=1-\alpha$. Here, we rely on the elementary inequality $s\ln s\leq s^\nfrac{5}{3}$ for $s>0$, the \GNI\ and the mass conservation \eqref{eq:mass-cons-n} to estimate
\begin{align*}
\intomega n_\epsi\ln n_\epsi\leq\big\|(n_\epsi+\epsi)^{1-\frac{\alpha}{2}}\big\|_{\Lo[\frac{5}{3(1-\frac{\alpha}{2})}]}^\frac{5}{3(1-\frac{\alpha}{2})}\leq C_1\big\|\nabla(n_\epsi+\epsi)^{1-\frac{\alpha}{2}}\big\|_{\Lo[2]}^\frac{5a}{3(1-\frac{\alpha}{2})}+C_1\quad\text{on }(0,\infty),
\end{align*}
with some $C_1>0$ and $a=\frac{12-6\alpha}{25-15\alpha}$. Since, in this case, $\alpha\leq 1$ we have $\frac{5a}{3(1-\frac{\alpha}{2})}\leq 2$ and hence (after an application of Young's inequality if necessary) there exists $C_2>0$ such that 
\begin{align*}
\frac{\intd}{\intd t}\intomega n_\epsi\ln n_\epsi+C_2\intomega n_\epsi\ln n_\epsi\leq \frac{S_0^2}{2m}\intomega |\nabla c_\epsi|^2+C_2\quad\text{on }(0,\infty).
\end{align*}
Due to Lemma \ref{lem:diffineq-lemma} and Lemma \ref{lem:bounds} this implies on one hand that there exists $C_3>0$ satisfying $\intomega n_\epsi\ln n_\epsi(\cdot,t)\leq C_3$ for all $t\geq0$ and on the other hand, upon returning to \eqref{eq:nab-est-alphaover2-eq3} and integrating with respect to time, that the asserted bound of the Lemma holds in light of the fact that $s\ln s\geq-\frac{1}{e}$ for all $s>0$.
\end{bew}

Now we can rely on standard reasoning to obtain the following.

\begin{lemma}\label{lem:time-reg-n-and-c}
Let $m\geq1$, $\alpha\geq0$ be such that $m+\alpha>\frac{4}{3}$ and assume that $n_0,c_0$ and $u_0$ comply with \eqref{IR} and that $S\in\CSp{2}{\bomega\times[0,\infty)^2;\R^{3\times 3}}$ fulfills \eqref{Sprop} with some $S_0>0$. For every $T>0$ there exists $C(T)>0$ such that for any $\epsi\in(0,1)$ the solution $(n_\epsi,c_\epsi,u_\epsi)$ of \eqref{approxprob} satisfies
\begin{align*}
\big\|\partial_t\big((n_{\epsi}+\epsi)^{m+\alpha-1}\big)\big\|_{\LSp{1}{(0,T);(W_0^{3,2}(\Omega))^*}}\leq C(T),
\end{align*}
and
\begin{align*}
\|c_{\epsi t}\|_{\LSp{1}{(0,T);(W_0^{3,2}(\Omega))^*}}\leq C(T).
\end{align*}
\end{lemma}

\begin{bew}
For fixed $T>0$ we find $C_1>0$ such that
\begin{align*}
\|\varphi\|_{\LSp{\infty}{(0,T);\W[1,\infty]}}\leq C_1\|\varphi\|_{\LSp{\infty}{(0,T);W_0^{3,2}(\Omega)}}\quad\text{for all }\varphi\in\LSp{\infty}{(0,T);W_0^{3,2}(\Omega)},
\end{align*}
in light of the continuous embedding of $\W[3,2]\hookrightarrow\W[1,\infty]$. Noting that $\LSp{\infty}{(0,T);W_0^{3,2}(\Omega)}$ is the dual space of $\LSp{1}{(0,T);\big(W_0^{3,2}(\Omega)\big)^*}$, we fix an arbitrary $\varphi\in\LSp{\infty}{(0,T);W_0^{3,2}(\Omega)}$ satisfying $\|\varphi\|_{\LSp{\infty}{(0,T);W_0^{3,2}(\Omega)}}\leq 1$ and make use of the first equation of \eqref{approxprob}, the Cauchy--Schwarz inequality and the bound \eqref{Sprop} to obtain
\begin{align*}
&\frac{1}{m+\alpha-1}\Big|\intomega\partial_t\big((n_\epsi+\epsi)^{m+\alpha-1}\big)\varphi\Big|\\
\leq\ &\frac{m|m+\alpha-2|C_1}{(m+\frac{\alpha}{2}-1)^2}\intomega\big|\nabla(n_\epsi+\epsi)^{m+\frac{\alpha}{2}-1}\big|^2\\
&\quad+\frac{mC_1}{m+\frac{\alpha}{2}-1}\Big(\intomega(n_\epsi+\epsi)^{2(m+\frac{\alpha}{2}-1)}\Big)^\frac{1}{2}\Big(\intomega\big|\nabla(n_\epsi+\epsi)^{m+\frac{\alpha}{2}-1}\big|^2\Big)^\frac{1}{2}\\
&\qquad+\frac{|m+\alpha-2|S_0C_1}{m+\frac{\alpha}{2}-1}\Big(\intomega\frac{n_\epsi^2(n_\epsi+\epsi)^{\alpha-2}}{(1+\epsi n_\epsi)^6(1+n_\epsi)^{2\alpha}}\big|\nabla(n_\epsi+\epsi)^{m+\frac{\alpha}{2}-1}\big|^2\Big)^\frac{1}{2}\Big(\intomega|\nabla c_\epsi|^2\Big)^\frac{1}{2}\\
&\quad\qquad+S_0C_1\Big(\intomega\frac{n_\epsi^2(n_\epsi+\epsi)^{2(m+\alpha-2)}}{(1+\epsi n_\epsi)^6(1+n_\epsi)^{2\alpha}}\Big)^\frac{1}{2}\Big(\intomega|\nabla c_\epsi|^2\Big)^\frac{1}{2}\\
&\qquad\qquad+\frac{C_1}{m+\alpha-1}\Big(\intomega|u_\epsi|^2\Big)^\frac{1}{2}\Big(\intomega\big|\nabla(n_\epsi+\epsi)^{m+\alpha-1}\big|^2\Big)^\frac{1}{2}\quad\text{on }(0,T)\text{ for all }\epsi\in(0,1).
\end{align*}
Since $\frac{n_\epsi^2(n_\epsi+\epsi)^{\alpha-2}}{(1+\epsi n_\epsi)^6(1+n_\epsi)^{2\alpha}}\leq \frac{(n_\epsi+\epsi)^\alpha}{(1+n_\epsi)^{2\alpha}}\leq 1$, multiple applications of the Young inequality and integration over $(0,T)$ entails the existence of $C_2>0$ such that
\begin{align*}
\intoT\Big|\intomega\partial_t\big((n_\epsi+\epsi)^{m+\alpha-1}\big)\varphi\Big|&\leq C_2\intoTomega\big|\nabla(n_\epsi+\epsi)^{m+\frac{\alpha}{2}-1}\big|^2+C_2\intoTomega\big|\nabla(n_\epsi+\epsi)^{m+\alpha-1}\big|^2\\&\qquad+C_2\intoTomega|\nabla c_\epsi|^2+C_2\intoTomega(n_\epsi+\epsi)^{2m+\alpha-2}+C_2\intoTomega|u_\epsi|^2+C_2
\end{align*}
holds for all $\epsi\in(0,1)$ and all $\varphi\in\LSp{\infty}{(0,T);W_0^{3,2}(\Omega)}$ with $\|\varphi\|_{\LSp{\infty}{(0,T);W_0^{3,2}(\Omega)}}\leq 1$. Because of $2m+\alpha-2<2(m+\alpha)-\frac{4}{3}$, a combination of Lemmata \ref{lem:bounds}, \ref{lem:st-bound-nepsi}, \ref{lem:epsi-ind-est-u} and \ref{lem:nab-est-alphaover2} now leads to the existence of $C_3(T)>0$ satisfying
\begin{align*}
\intoT\Big|\intomega\partial_t\big((n_\epsi+\epsi)^{m+\alpha-1}\big)\varphi\Big|\leq C_3(T)\ \text{for all }\varphi\in\LSp{\infty}{(0,T);W_0^{2,3}(\Omega)}\text{ with }\|\varphi\|_{\LSp{\infty}{(0,T);W_0^{3,2}(\Omega)}}\leq 1.
\end{align*}
For the second part of the Lemma we follow a follow complementary reasoning for the second equation. For fixed $\varphi$ as before we obtain $C_4>0$ such that
\begin{align*}
\Big|\intomega c_{\epsi t}\varphi\Big|\leq C_1\intomega|\nabla c_\epsi|^2+C_1\intomega c_\epsi+C_1\intomega n_\epsi+\frac{C_1}{2}\intomega |u_\epsi|^2+\frac{C_1}{2}\intomega c_\epsi^2+C_4
\end{align*}
is valid in $(0,T)$ for all $\epsi\in(0,1)$. Hence, we can conclude the proof upon integration over $(0,T)$ in light of the bounds featured in Lemmata \ref{lem:glob_ex}, \ref{lem:bounds} and \ref{lem:epsi-ind-est-u}.
\end{bew}

Enhancing arguments akin to those present in the previous proof by known results for the Yosida approximation and the Stokes operator, a similar result can be established for the third solution component.

\begin{lemma}\label{lem:time-reg-u}
Let $m\geq1$, $\alpha\geq0$ be such that $m+\alpha>\frac{4}{3}$ and suppose that $n_0,c_0$ and $u_0$ fulfill \eqref{IR} and that $S\in\CSp{2}{\bomega\times[0,\infty)^2;\R^{3\times 3}}$ satisfies \eqref{Sprop} with some $S_0>0$. For every $T>0$ there exists $C(T)>0$ such that for any $\epsi\in(0,1)$ the solution $(n_\epsi,c_\epsi,u_\epsi)$ of \eqref{approxprob} satisfies
\begin{align}\label{eq:time-reg-u}
\int_0^T\|u_{\epsi t}\|^{\frac43}_{(W_{0,\sigma}^{1,2}(\Omega))^*}\leq C(T).
\end{align}
\end{lemma}

\begin{bew}
In light of \eqref{eq:st-bound-specialcases} from Lemma \ref{lem:st-bound-nepsi} there is $C_1>0$ such that $\|n_\epsi(\cdot,t)\|_{\Lo[\frac65]}^\frac43\leq C_1$ for all $t>0$ and hence we can follow the proof of \cite[Lemma 5.5]{Wang17-globweak-ksns}, where the related system with linear diffusion was discussed, to conclude the desired bound. Let us state a brief outline of the steps involved. We multiply the third equation of \eqref{approxprob} with a fixed $\psi\in C_0^\infty(\Omega)$ satisfying $\nabla\cdot\psi\equiv0$ throughout $\Omega$ and employ Hölder's inequality to obtain
\begin{align*}
\Big|\intomega u_{\epsi t}\cdot\psi\Big|\leq\|\nabla u_\epsi\|_{\Lo[2]}\|\nabla\psi\|_{\Lo[2]}+\|Y_\epsi u_\epsi\|_{\Lo[6]}&\|u_\epsi\|_{\Lo[3]}\|\nabla u_\epsi\|_{\Lo[2]}\\&+\|\nabla\phi\|_{\Lo[\infty]}\|n_\epsi\|_{\Lo[\frac{6}{5}]}\|\psi\|_{\Lo[6]}
\end{align*}
on $(0,\infty)$ for all $\epsi\in(0,1)$. Next, we make use of known facts for the Yoshida approximation and the Stokes operator, the embedding $W_{0,\sigma}^{1,2}(\Omega)\hookrightarrow\Lo[6]$ and the \GNI\ to obtain $C_2>0$ such that
\begin{align*}
\|Y_\epsi u_\epsi(\cdot,t)\|_{\Lo[6]}\leq\|\nabla u_\epsi(\cdot,t)\|_{\Lo[2]},\quad\text{and}\quad \|u_\epsi(\cdot,t)\|_{\Lo[3]}^\frac{4}{3}\leq C_2\|\nabla u_\epsi(\cdot,t)\|_{\Lo[2]}^\frac{2}{3}\|u_\epsi(\cdot,t)\|_{\Lo[2]}^\frac{2}{3}
\end{align*}
for all $t>0$ and all $\epsi\in(0,1)$. Combining the estimates above with Young's inequality shows that with some $C_3>0$ we have
\begin{align*}
\intoT\|u_{\epsi t}\|_{(W_{0,\sigma}^{1,2}(\Omega))^*}^\frac{4}{3}\leq C_3\intoT\|\nabla u_\epsi\|_{\Lo[2]}^\frac{4}{3}+C_3\intoT\|\nabla u_\epsi\|^2_{\Lo[2]}\|u_\epsi\|_{\Lo[2]}^\frac{2}{3}+C_3T
\end{align*}
for all $T>0$ and all $\epsi\in(0,1)$, completing the proof in terms of Lemma \ref{lem:epsi-ind-est-u}. 
\end{bew}

\setcounter{equation}{0} 
\section{Limit functions and their regularity properties}\label{sec5:convergence}
The uniform bounds prepared in the previous section enabled us to derive the existence limit functions $n,c,u$ satisfying the regularity conditions imposed by Definition \ref{def:very_weak_sol}. In addition, the precompactness properties contained in the Lemmata of the previous section will enable us to pass to obtain convergence properties suitable for passing to the limit in most of the integrals making up the solution concepts discussed in Section \ref{sec2:soldef}. In contrast to the scalar sensitivity case discussed in \cite{TB2017_nonlindiff} and the linear diffusion case discussed in \cite{Wang17-globweak-ksns} the very weak solution concept features terms combining $n_\epsi+1$ and $n_\epsi+\epsi$ in a slightly more varied way, necessitating the preparation of additional convergence properties.

\begin{lemma}\label{lem:convergences}
Let $m\geq1$, $\alpha\geq0$ be such that $m+\alpha>\frac{4}{3}$ and suppose that $n_0,c_0,u_0$ comply with \eqref{IR} and assume that $S\in\CSp{2}{\bomega\times[0,\infty)^2;\R^{3\times 3}}$ fulfills \eqref{Sprop} with some $S_0>0$.  Then there exist a sequence $(\epsi_j)_{j\in\N}\subset(0,1)$ with $\epsi_j\searrow0$ as $j\to\infty$ and functions
\begin{align*}
n&\in\LSploc{2(m+\alpha)-\frac{4}{3}}{\bomega\times[0,\infty)}\quad\text{with}\quad\nabla n^{m+\alpha-1}\in\LSploc{2}{\bomega\times[0,\infty)},\\
c&\in\LSploc{2}{[0,\infty);\W[1,2]},\\
u&\in\LSploc{2}{[0,\infty);W^{1,2}_{0,\sigma}(\Omega)},
\end{align*}
such that the solutions $(n_\epsi,c_\epsi,u_\epsi)$ of \eqref{approxprob} satisfy
\begin{alignat}{2}
(n_\epsi+\epsi)^{m+\alpha-1}&\to n^{m+\alpha-1}&&\text{in }\LSploc{2}{\bomega\times[0,\infty)}\text{ and a.e. in }\Omega\times(0,\infty),\label{eq:conv-n-m+a-1-ae}\\
\nabla (n_\epsi+\epsi)^{m+\alpha-1}&\wto \nabla n^{m+\alpha-1}\ &&\text{in }\LSploc{2}{\bomega\times[0,\infty)},\label{eq:conv-nab-n-m+a-1}\\
n_\epsi+\epsi&\wto n\quad\qquad&&\text{in }\LSploc{2(m+\alpha)-\frac{4}{3}}{\bomega\times[0,\infty)},\label{eq:conv-nw}\\
n_\epsi+\epsi\to n\quad\text{and}\quad n_\epsi&\to n&&\text{in }\LSploc{p}{\bomega\times[0,\infty)}\text{ for any }p\in[1,2(m+\alpha)-\tfrac43),\label{eq:conv-n-strong}\\
c_\epsi&\to c&&\text{in }\LSploc{2}{\bomega\times[0,\infty)}\text{ and a.e in }\Omega\times(0,\infty),\label{eq:conv-c}\\
\nabla c_\epsi&\wto \nabla c&&\text{in }\LSploc{2}{\bomega\times[0,\infty)},\label{eq:conv-nab-c}
\intertext{as well as}
u_\epsi&\to u&&\text{in }\LSploc{2}{\bomega\times[0,\infty)}\text{ and a.e. in }\Omega\times(0,\infty),\label{eq:conv-u-l2}\\
\nabla u_\epsi&\wto \nabla u&&\text{in }\LSploc{2}{\bomega\times[0,\infty)},\label{eq:conv-nab-u}\\
Y_\epsi u_\epsi&\to u&&\text{in }\LSploc{2}{\bomega\times[0,\infty)}\label{eq:conv-Yu}\\
n_\epsi u_\epsi&\to nu&&\text{in }\LSploc{1}{\bomega\times[0,\infty)}\label{eq:conv-nu-l1}
\end{alignat}
as $\epsi=\epsi_j\searrow0$, and such that $n\geq0$, $c\geq0$ a.e. in $\Omega\times(0,\infty)$. If, moreover, $m+\alpha\in(\tfrac43,2]$, then there exists a further subsequence $(\epsi_{j_k})_{k\in\N}\subset(0,1)$ such that $(n_\epsi,c_\epsi,u_\epsi)$ also satisfy
\begin{alignat}{2}
(n_\epsi+1)^{m+\alpha-1}&\to (n+1)^{m+\alpha-1}&&\text{in }\LSploc{2}{\bomega\times[0,\infty)},\label{eq:conv-n+1-m+a-1-ae}\\
(n_\epsi+1)^{m+2\alpha-1}&\to (n+1)^{m+2\alpha-1}&&\text{in }\LSploc{1}{\bomega\times[0,\infty)},\label{eq:conv-n+1-m+2alpha-1}\\
\nabla (n_\epsi+1)^{m+\alpha-1}&\wto \nabla (n+1)^{m+\alpha-1}&&\text{in }\LSploc{2}{\bomega\times[0,\infty)},\label{eq:conv-nab-n+1-m-1}\\
(n_\epsi+1)^{\frac{m+2\alpha-3}{2}}(n_\epsi+\epsi)^{\frac{m-1}{2}}\nabla n_\epsi&\wto (n+1)^{\frac{m+2\alpha-3}{2}}n^{\frac{m-1}{2}} \nabla n\quad&&\text{in }\LSploc{2}{\bomega\times[0,\infty)},\label{eq:conv-nab-n+1+eps-m-1}\\
(n_\epsi+1)^\alpha(n_\epsi+\epsi)^{m-1}&\to(n+1)^\alpha n^{m-1}&&\text{in }\LSploc{2}{\bomega\times[0,\infty)},\label{eq:conv-n+1+eps-alpha-m-1}\\
(n_\epsi+1)^{m+2\alpha-1}u_\epsi&\to(n+1)^{m+2\alpha-1}u&&\text{in }\LSploc{1}{\bomega\times[0,\infty)},\label{eq:conv-n+1-m+2alpha-1-u-l1}
\end{alignat}
as $\epsi=\epsi_{j_k}\searrow0$.
\end{lemma}

\begin{bew}
Noticing that $2(m+\alpha-1)<2(m+\alpha)-\frac{4}{3}$, we find that by combining Lemmata \ref{lem:bounds}, \ref{lem:st-bound-nepsi} and \ref{lem:time-reg-n-and-c} with the Aubin-Lions lemma (\cite[Corollary 8.4]{Sim87}) 
\begin{align*}
\big\{(n_\epsi+\epsi)^{m+\alpha-1}\big\}_{\epsi\in(0,1)}\quad\text{is relative compact in }\LSploc{2}{\bomega\times[0,\infty)}
\end{align*}
and that hence there exists a sequence $\epsi_j\searrow0$  such that \eqref{eq:conv-n-m+a-1-ae} holds. Extracting an additional subsequence (still denoted by $\epsi_j$) we conclude from the spatio-temporal bounds featured in Lemma \ref{lem:bounds} and Lemma \ref{lem:st-bound-nepsi} that \eqref{eq:conv-nab-n-m+a-1} and \eqref{eq:conv-nw} hold as well. In light of Lemma \ref{lem:st-bound-nepsi} $\{(n_{\epsi_j}+\epsi_j)^p\}_{j\in\N}$ is equi-integrable for any $p<2(m+\alpha)-\frac43$, and thus we can rely on the a.e. convergence of $n_\epsi+\epsi$ entailed by \eqref{eq:conv-n-m+a-1-ae} and the Vitali convergence theorem to obtain the first part of \eqref{eq:conv-n-strong}, with the second part then being an immediate consequence of the uniform convergence of $\epsi_j$ to zero. Along similar lines the Lemmata \ref{lem:bounds} and \ref{lem:time-reg-n-and-c} together with the Aubin-Lions lemma imply that upon extraction of another subsequence also \eqref{eq:conv-c} and \eqref{eq:conv-nab-c} hold. Moreover, applying these arguments once more for the third component of the approximate solution while relying on Lemmata \ref{lem:epsi-ind-est-u} and \ref{lem:time-reg-u} proves \eqref{eq:conv-u-l2} and \eqref{eq:conv-nab-u}, whereas \eqref{eq:conv-Yu} is a consequence of the dominated convergence theorem and the boundedness of $\|u_\epsi\|_{\LSp{2}{\Omega\times(0,\infty)}}^2$ (see e.g. \cite[Lemma 4.1]{win_globweak3d-AHPN16}). The strong convergence property of the mixed term $n_\epsi u_\epsi$ in \eqref{eq:conv-nu-l1} can be concluded by combining the a.e. convergences contained in \eqref{eq:conv-n-m+a-1-ae} and \eqref{eq:conv-u-l2} with the equi-integrability of $\{|n_{\epsi_j} u_{\epsi_j}|^r\}_{j\in\N}$ for some $r>1$ implied by Lemma \ref{lem:nepsi-uepsi-bounds} and Vitali's convergence theorem. The assertions for the special case of $m+\alpha\in(\frac43,2)$ follow from identical reasoning in light of Corollary \ref{cor:1} and Lemma \ref{lem:nepsi-uepsi-bounds}. To be precise, we can conclude (upon extraction of another non-relabeled subsequence) \eqref{eq:conv-nab-n+1-m-1} and \eqref{eq:conv-nab-n+1+eps-m-1} from \eqref{eq:cor1}. The properties \eqref{eq:conv-n+1-m+a-1-ae}, \eqref{eq:conv-n+1-m+2alpha-1} and \eqref{eq:conv-n+1+eps-alpha-m-1} are a consequence of \eqref{eq:cor2}, Vitali's convergence theorem and the fact that $m+2\alpha-1\leq 2(m+\alpha-1)$, and finally, combining Lemma \ref{lem:nepsi-uepsi-bounds} with Vitali's theorem one last time shows \eqref{eq:conv-n+1-m+2alpha-1-u-l1}.
\end{bew}

\setcounter{equation}{0} 
\section{Solution properties of the limit functions}\label{sec6:sol-prop}

\subsection{Weak solution properties of \texorpdfstring{$c$ and $u$}{c and u}}\label{sec61:weak-sol-cu}
Reyling on the convergence properties prepared in Lemma \ref{lem:convergences}, we can check in a straightforward manner that the limit objects $c$ and $u$ are weak solutions of their corresponding equations in \eqref{CTndS}.

\begin{lemma}\label{lem:sol-prop-c-u}
Let $m\geq1$, $\alpha\geq0$ be such that $m+\alpha>\frac{4}{3}$, assume that $n_0,c_0$ and $u_0$ comply with \eqref{IR} and suppose that $S\in\CSp{2}{\bomega\times[0,\infty)^2;\R^{3\times 3}}$ satisfies \eqref{Sprop} with some $S_0>0$. Furthermore, let $n,c,u$ denote the limit functions provided by Lemma \ref{lem:convergences}. Then 
\begin{align}\label{eq:mass-eq}
\intomega n(\cdot,t)= \intomega n_0\quad\text{for a.e. }t>0,
\end{align}
and $c$ and $u$ satisfy the weak solution properties \eqref{eq:very_weak_sol1} and \eqref{eq:very_weak_sol2}, respectively, of Definition \ref{def:very_weak_sol}.
\end{lemma}

\begin{bew}
The equality in \eqref{eq:mass-eq} for almost every $t>0$ is a direct result of the mass conservation \eqref{eq:mass-cons-n} from Lemma \ref{lem:glob_ex} and \eqref{eq:conv-n-strong}. To verify that $c$ solves its corresponding equation in the weak sense, we multiply the second equation of \eqref{approxprob} by an arbitrary test function $\varphi\in\LSp{\infty}{\Omega\times(0,\infty)}\cap\LSp{2}{(0,\infty);\W[1,2]}$ with compact support in $\bomega\times[0,\infty)$ and $\varphi_t\in\LSp{2}{\Omega\times(0,\infty)}$ to find that
\begin{align*}
-\intinfomega c_\epsi\varphi_t-\intomega c_0\varphi(\cdot,0)=-\intinfomega\nabla c_\epsi\cdot\nabla \varphi-\intinfomega n_\epsi\varphi+\intinfomega c_\epsi(u_\epsi\cdot\nabla\varphi)
\end{align*}
holds for all $\epsi\in(0,1)$. In consideration of \eqref{eq:conv-c}, \eqref{eq:conv-nab-c}, \eqref{eq:conv-n-strong} and \eqref{eq:conv-u-l2} we may pass to the limit in each of the integrals and conclude that \eqref{eq:very_weak_sol1} holds and that hence $c$ solves its equation in the weak sense. In a similar fashion, we test the third equation of \eqref{approxprob} by an arbitrary $\psi\in C_0^\infty\left(\Omega\times[0,\infty)\right)$ satisfying $\nabla\cdot\psi\equiv0$ in $\Omega\times(0,\infty)$ to obtain
\begin{align*}
-\intinfomega u_\epsi\psi_t-\intomega u_0\psi(\cdot,0)=-\intinfomega\nabla u_\epsi\cdot\nabla\psi+\intinfomega(Y_\epsi u_\epsi\otimes u_\epsi)\cdot\nabla \psi+\intinfomega n_\epsi(\nabla\phi\cdot\psi)
\end{align*}
for all $\epsi\in(0,1)$. Recalling \eqref{eq:conv-u-l2}, \eqref{eq:conv-nab-u}, \eqref{eq:conv-Yu}, as well as \eqref{eq:conv-n-strong} and \eqref{phidef} we can take $\epsi\searrow0$ in all the integrals and find that $u$ satisfies \eqref{eq:very_weak_sol2}.
\end{bew}

\subsection{Weak solution property of \texorpdfstring{$n$ for $m+2\alpha>\frac53$}{n for m+2alpha>5/3}}\label{sec62:weak-sol-n}
The currently known compactness properties do not allow us to take $\epsi\searrow0$ in some of the integrals appearing in the equation for $n_\epsi$ corresponding to \eqref{eq:weak-sol-n} of the weak solution concept in Definition \ref{def:weak-solution}. However, imposing the additional condition $m+2\alpha>\frac{5}{3}$ we can obtain supplementary convergence properties to the ones in Lemma \ref{lem:convergences}, which will allow us to pass to the limit in these crucial integrals.

\begin{lemma}\label{lem:weak-sol-5over3}
Let $m\geq1$, $\alpha\geq0$ be such that $m+2\alpha>\frac{5}{3}$, suppose that $n_0,c_0$ and $u_0$ comply with \eqref{IR}, and suppose that $S\in\CSp{2}{\bomega\times[0,\infty)^2;\R^{3\times 3}}$ satisfies \eqref{Sprop} with some $S_0>0$. Furthermore, let $n,c,u$ denote the limit functions obtained in Lemma \ref{lem:convergences}. Then $n\in\LSploc{2}{\bomega\times[0,\infty)}$ and for any $\varphi\in C_0^\infty\big(\bomega\times[0,\infty)\big)$ the weak solution property \eqref{eq:weak-sol-n} is satisfied.
\end{lemma}

\begin{bew}
Multiplying the first equation of \eqref{approxprob} by $\varphi\in C_0^\infty\left(\bomega\times[0,\infty)\right)$ and integrating by parts, we find that
\begin{align}\label{eq:weak-form-nepsi}
-\intinfomega n_\epsi\varphi_t-\intomega n_0\varphi(\cdot,0)=-&\frac{m}{m+\alpha-1}\intinfomega(n_\epsi+\epsi)^{1-\alpha}\big(\nabla (n_\epsi+\epsi)^{m+\alpha-1}\cdot\nabla\varphi\big)\\&+\intinfomega \frac{n_\epsi}{(1+\epsi n_\epsi)^3}\big(S_\epsi(x,n_\epsi,c_\epsi)\nabla c_\epsi\cdot\nabla \varphi\big)+\intinfomega n_\epsi(u_\epsi\cdot\nabla\varphi)\nonumber
\end{align}
holds for all $\epsi\in(0,1)$, where we used $(n_\epsi+\epsi)^{m-1}\nabla n_\epsi=\frac{(n_\epsi+\epsi)^{1-\alpha}}{m+\alpha-1}\nabla (n_\epsi+\epsi)^{m+\alpha-1}$. In light of \eqref{eq:conv-n-strong} we see that
\begin{align*}
-\intinfomega n_\epsi\varphi_t\to-\intinfomega n\varphi_t\quad\text{as}\ \epsi=\epsi_j\searrow0.
\end{align*}
Moreover, since $m+2\alpha>\frac{5}{3}$, we have $2(1-\alpha)<2(m+\alpha)-\frac43$, so that \eqref{eq:conv-n-strong} implies that
\begin{align*}
(n_\epsi+\epsi)^{1-\alpha}\to n^{1-\alpha}\quad\text{in}\quad\LSploc{2}{\bomega\times[0,\infty)}\quad \text{as}\quad \epsi=\epsi_j\searrow0,
\end{align*}
which together with \eqref{eq:conv-nab-n-m+a-1} shows
\begin{align*}
-\frac{m}{m+\alpha-1}\intinfomega (n_\epsi&+\epsi)^{1-\alpha}\big(\nabla (n_\epsi+\epsi)^{m+\alpha-1}\cdot\nabla \varphi\big)\\&\to -\frac{m}{m+\alpha-1}\intinfomega n^{1-\alpha}\big(\nabla n^{m+\alpha-1}\cdot\nabla \varphi\big)=- m \intinfomega n^{m-1}(\nabla n\cdot\nabla\varphi)
\end{align*}	
as $\epsi=\epsi_j\searrow0$.  Additionally, since $2(1-\alpha)<2(m+\alpha)-\frac43$, we can fix $2<s<\frac{2(m+\alpha)-\frac43}{1-\alpha}$ and find that
\begin{align*}
\int_t^{t+1}\!\intomega\Big|\frac{n_\epsi S_\epsi(x,n_\epsi,c_\epsi)}{(1+\epsi n_\epsi)^3}\Big|^s \leq \int_t^{t+1}\!\intomega\frac{S_0^s n_\epsi^s}{(1+n_\epsi)^{s\alpha}}\leq\begin{cases}
S_0^s|\Omega|,\quad&\text{if }\alpha\geq1\\S_0^2
\int_t^{t+1}\!\intomega(n_\epsi+\epsi)^{s(1-\alpha)},&\text{if }\alpha\in[0,1)
\end{cases}
\end{align*}
holds on $(0,\infty)$. Making use of the fact that $s(1-\alpha)<2(m+\alpha)-\frac43$ and Lemma \ref{lem:st-bound-nepsi} we thus obtain that $\big\{n_{\epsi_j}^2 S_{\epsi_j}(x,n_{\epsi_j},c_{\epsi_j})^2(1+{\epsi_j} n_{\epsi_j})^{-6}\big\}_{j\in\N}$ is equi-integrable, which together with the a.e. convergences of $S_\epsi\to S$ and $\frac{n_\epsi}{(1+\epsi n_\epsi)^3}\to n$ in $\Omega\times(0,\infty)$ and Vitali's theorem shows that
\begin{align*}
\frac{n_\epsi S_\epsi(x,n_\epsi,c_\epsi)}{(1+\epsi n_\epsi)^3}\to n S(x,n,c)\quad\text{in}\quad\LSploc{2}{\bomega\times[0,\infty)}
\end{align*}
as $\epsi=\epsi_j\searrow0$. Merging this convergence property with \eqref{eq:conv-nab-c} we obtain that
\begin{align*}
\intinfomega\frac{n_\epsi}{(1+\epsi n_\epsi)^3}\big(S_\epsi(x,n_\epsi,c_\epsi)\nabla c_\epsi\cdot\nabla \varphi\big)\to\intinfomega n\big(S(x,n,c)\nabla c\cdot\nabla\varphi\big)\quad\text{as}\quad \epsi=\epsi_j\searrow0.
\end{align*}
Finally, relying on \eqref{eq:conv-nu-l1} we see that 
\begin{align*}
\intinfomega n_\epsi(u_\epsi\cdot\nabla\varphi)\to\intinfomega n(u\cdot\nabla\varphi)\quad\text{as}\quad \epsi=\epsi_j\searrow0.
\end{align*}
In conclusion, we may pass to the limit in each of the integrals in \eqref{eq:weak-form-nepsi} and find that \eqref{eq:weak-sol-n} holds.
\end{bew}

Amalgamating the previous results finalizes the proof of Theorem \ref{theo:1}.

\begin{proof}[\textbf{Proof of Theorem \ref{theo:1}}:]
The proof is immediate after combination of Lemmata \ref{lem:sol-prop-c-u} and \ref{lem:weak-sol-5over3} with the regularity information on $n,c$ and $u$ presented in Lemma \ref{lem:convergences}.
\end{proof}

\subsection{Very weak solution property of \texorpdfstring{$n$}{n} in the case of \texorpdfstring{$m+\alpha>\frac{4}{3}$}{m+alpha>4/3}}\label{sec63:very-weak-sol}
Under the weaker assumption that only $m+\alpha>\frac{4}{3}$ the obtained limit function $n$ does not appear to be regular enough to conclude that the integral $\intinfomega n^{m-1}\nabla n\cdot\varphi$, appearing in \eqref{eq:weak-sol-n}, is finite. Weakening the solution concept appears to be the only way to compensate the missing regularity information, which is why we will only consider global very weak solutions as defined in Definition \ref{def:very_weak_sol} for the parameter range of $m\geq1$ and $\alpha\geq0$ between $m+\alpha>\frac{4}{3}$ and $m+2\alpha\leq \frac{5}{3}$. Working under these weaker hypothesis, however, the weak convergence statement for $\nabla c_\epsi$ is insufficient to pass to the limit in the integral containing both gradient terms. Therefore, we will have to attain a strong convergence result for $\nabla c_\epsi$ which we prepare with the following Lemma from \cite{Wang17-globweak-ksns}.

\begin{lemma}\label{lem:c2-ineq}
Let $m\geq1$, $\alpha\geq0$ be such that $m+\alpha>\frac{4}{3}$ and assume that $n_0,c_0$ and $u_0$ comply with \eqref{IR} and suppose that $S\in\CSp{2}{\bomega\times[0,\infty)^2;\R^{3\times 3}}$ satisfies \eqref{Sprop} with some $S_0>0$. Then there exists a null set $N\subset(0,\infty)$ such that the functions $n,c$ and $u$ obtained in Lemma \ref{lem:convergences} satisfy
\begin{align}\label{eq:c2-ineq}
\frac{1}{2}\intomega c^2(\cdot,T)-\frac{1}{2}\intomega c_0^2+\intoTomega|\nabla c|^2\geq -\intoTomega c^2+\intoTomega nc\quad\text{for all }T\in(0,\infty)\setminus N.
\end{align}
\end{lemma}

\begin{bew}
This is precisely \cite[Lemma 7.1]{Wang17-globweak-ksns}. The same lemma has also been used in the setting with scalar sensitivity in \cite[Lemma 6.3]{TB2017_nonlindiff}. We will refrain from repeating the rather technical argumentation concerning Steklov averages underlying the proof and refer the reader to \cite{Wang17-globweak-ksns} for details.
\end{bew}

Relying on the spatio-temporal estimates of Section \ref{sec4:reg-est} and the inequality above we can now pass to another subsequence along which $\nabla c_\epsi\to\nabla c$ in $\LSp{2}{\Omega\times(0,T)}$ holds as $\epsi\searrow0$. Similar reasoning has been employed in e.g. \cite[Lemma 4.4]{Winkler18_stokesrot} and \cite[Lemma 7.2]{Wang17-globweak-ksns}.

\begin{lemma}\label{lem:strong-conv-nabc-l2}
Let $m\geq1$, $\alpha\geq0$ be such that $m+\alpha>\frac{4}{3}$ and assume that $n_0,c_0$ and $u_0$ comply with \eqref{IR} and suppose that $S\in\CSp{2}{\bomega\times[0,\infty)^2;\R^{3\times 3}}$ satisfies \eqref{Sprop} with some $S_0>0$. Furthermore, denote by $(\epsi_j)_{j\in\N}$ and $n,c,u$ the sequence and limit functions provided by Lemma \ref{lem:convergences}. Then there exist a subsequence $(\epsi_{j_k})_{k\in\N}$ and a null set $N\subset(0,\infty)$ such that for each $T\in(0,\infty)\setminus N$ the classical solution $(n_\epsi,c_\epsi,u_\epsi)$ of \eqref{approxprob} satisfies
\begin{align*}
\nabla c_\epsi\to\nabla c\quad\text{in }\LSp{2}{\Omega\times(0,T)}\ \text{as}\ \epsi=\epsi_{j_k}\searrow0.
\end{align*}
\end{lemma}

\begin{bew}
With $r\in(1,2)$ as given by Lemma \ref{lem:st-bound-nepsi} we note that, due to the bounds presented in Lemmata \ref{lem:bounds} and \ref{lem:st-bound-nepsi}, the nonnegativity of $n_\epsi$ and the Hölder and Young inequalities we have $C>0$ satisfying
\begin{align*}
\int_t^{t+1}\!\!\intomega|n_\epsi c_\epsi|^r&\leq\int_t^{t+1}\!\!\|n_\epsi+\epsi\|_{\Lo[\frac{6r}{6-r}]}^r\|c_\epsi\|_{\Lo[6]}^r
\leq\frac{2-r}{2}\int_t^{t+1}\!\!\|n_\epsi+\epsi\|_{\Lo[\frac{6r}{6-r}]}^\frac{2r}{2-r}+\frac{r}{2}\int_t^{t+1}\!\!\|c_\epsi\|_{\Lo[6]}^2\leq C
\end{align*}
for all $t>0$ and all $\epsi\in(0,1)$. Since $r>1$, we can combine the a.e. convergence of $n_\epsi c_\epsi\to nc$ in $\Omega\times(0,\infty)$ as $\epsi=\epsi_j\searrow0$, as implied by Lemma \ref{lem:convergences}, with Vitali's convergence theorem, to find that for all $T>0$
\begin{align*}
\intoTomega n_\epsi c_\epsi\to\intoTomega nc\quad\text{as }\epsi=\epsi_j\searrow0.
\end{align*}
Denoting by $N_1\subset(0,\infty)$ the null set given by Lemma \ref{lem:c2-ineq} we see that by Lemma \ref{lem:convergences} there exists another null set $N_2\supset N_1$ and a subsequence $(\epsi_{j_k})_{k\in\N}$ such that
\begin{align*}
\intomega c_\epsi^2(\cdot,T)\to\intomega c^2(\cdot,T)\quad\text{for all }T\in(0,\infty)\setminus N_2\quad\text{as }\epsi=\epsi_{j_k}\searrow0.
\end{align*}
Hence, for any such $T\in(0,\infty)\setminus N_2$, by testing the second equation of \eqref{approxprob} by $c_\epsi$ and making use of Lemmata \ref{lem:c2-ineq} and \ref{lem:convergences} we obtain
\begin{align*}
\intoTomega|\nabla c|^2&\geq-\frac{1}{2}\intomega c^2(\cdot,T)+\frac{1}{2}\intomega c_0^2-\intoTomega c^2+\intoTomega nc\\
&=\lim_{\epsi_{j_k}\searrow0}\Big(-\frac{1}{2}\intomega c^2_{\epsi_{j_k}}+\frac{1}{2}\intomega c_0^2-\intoTomega c_{\epsi_{j_k}}^2+\intoTomega n_{\epsi_{j_k}}c_{\epsi_{j_k}}\Big)=\lim_{\epsi_{j_k}\searrow0}\intoTomega|\nabla c_{\epsi_{j_k}}|^2,
\end{align*}
which together with the fact that the norm in $\LSp{2}{\Omega\times(0,T)}$ is weakly lower semicontinuous and the weak convergence property in \eqref{eq:conv-nab-c} implies that actually $\nabla c_\epsi\to\nabla c$ in $\LSp{2}{\Omega\times(0,T)}$ as $\epsi=\epsi_{j_k}\searrow0$.
\end{bew}

Finally, as a last step before proving Theorem \ref{theo:2}, we can verify the $\Phi$-supersolution property of Definition \ref{def:weak_super_sol} for the choice of $\Phi(s)=(s+1)^{m+2\alpha-1}$ whenever $m\geq1$ and $\alpha\geq0$ satisfy $m+\alpha>\frac{4}{3}$ and $m+2\alpha<2$. The restriction $m+2\alpha<2$, however, is of no consequence for our Theorem, since for $m\geq1$ and $\alpha\geq0$ with $m+\alpha>\frac43$ and $m+2\alpha\geq2$ the existence of a global very weak solution is already established by Theorem \ref{theo:1} in light of the fact that every weak solution is also a very weak solution.

\begin{lemma}\label{lem:sol-prop-n}
Let $m\geq1$, $\alpha\geq0$ satisfy $m+\alpha>\frac43$ and $m+2\alpha<2$. Assume that $n_0,c_0,u_0$ comply with \eqref{IR} and suppose that $S\in\CSp{2}{\bomega\times[0,\infty)^2;\R^{3\times 3}}$ fulfills \eqref{Sprop} with some $S_0>0$. Moreover, denote by $n,c,u$ the limit functions provided by Lemma \ref{lem:convergences} and let $\Phi(s):=(s+1)^{m+2\alpha-1}$ for $s\geq0$. Then $n$ is a global $\Phi$--supersolution of \eqref{CTndS} in the sense of Definition \ref{def:weak_super_sol}.
\end{lemma}

\begin{bew}
Because of $m+2\alpha<2$ and $\alpha\geq0$ we clearly have $m+\alpha\leq 2$ and we may hence draw on the special case convergences discussed in Lemma \ref{lem:convergences}, i.e. \eqref{eq:conv-n+1-m+a-1-ae}--\eqref{eq:conv-n+1-m+2alpha-1-u-l1}. With $\Phi(s):=(s+1)^{m+2\alpha-1}$ for $s\geq0$, we find that the regularity requirements demanded by Definition \ref{def:weak_super_sol} were already obtained in Lemma \ref{lem:convergences}. In particular, we find that the conditions concerning $n$ contained in \eqref{eq:very_weak_supersol_regularity} are implied by \eqref{eq:conv-n+1-m+2alpha-1}, \eqref{eq:conv-nab-n+1+eps-m-1}, \eqref{eq:conv-nab-n+1-m-1} together with \eqref{eq:conv-n+1+eps-alpha-m-1}, \eqref{eq:conv-n+1-m+2alpha-1-u-l1}, \eqref{eq:conv-n+1-m+a-1-ae} and \eqref{eq:conv-nab-n+1-m-1}, respectively, where we also used the fact that $\frac{n(n+1)^{\alpha-1}}{(1+n)^\alpha}\in\LSploc{\infty}{\bomega\times[0,\infty)}$. What remains is the verification of \eqref{eq:very_weak_supersol}. We pick an arbitrary nonnegative test function $\varphi\in C_0^\infty\big(\bomega\times[0,\infty)\big)$ satisfying $\frac{\partial\varphi}{\partial\nu}=0$ on $\romega\times(0,\infty)$ and then fix $T>0$ such that $\varphi\equiv0$ in $\Omega\times(T,\infty)$. Keeping in mind that $m+2\alpha<2$, we multiply the first equation of \eqref{approxprob} with $(m+2\alpha-1)(n_\epsi+1)^{m+2\alpha-2}\varphi$, integrate by parts and rewrite some terms to obtain that
\begin{align}\label{eq:sol-prop-n-proof-eq1}
&-\intoTomega(n_\epsi+1)^{m+2\alpha-1}\varphi_t-\intomega(n_0+1)^{m+2\alpha-1}\varphi(\cdot,0)\nonumber\\
=&\ m(m+2\alpha-1)(2-(m+2\alpha))\intoTomega\big|(n_\epsi+1)^{\frac{m+2\alpha-3}{2}}(n_\epsi+\epsi)^{\frac{m-1}{2}}\nabla n_\epsi\big|^2\varphi\nonumber\\
&\quad  -\frac{m(m+2\alpha-1)}{m+\alpha-1}\intoTomega(n_\epsi+1)^{\alpha}(
n_\epsi+\epsi)^{m-1}\big(\nabla (n_\epsi+1)^{m+\alpha-1}\cdot\nabla\varphi\big)\\
&\qquad -\frac{(m+2\alpha-1)(2-(m+2\alpha))}{m+\alpha-1}\intoTomega\frac{(n_\epsi+1)^{\alpha-1}n_\epsi}{(1+\epsi n_\epsi)^3}\big(\nabla(n_\epsi+1)^{m+\alpha-1}\cdot S_\epsi(\cdot,n_\epsi,c_\epsi)\nabla c_\epsi\big)\varphi\nonumber\\
&\quad\qquad+(m+2\alpha-1)\intoTomega(n_\epsi+1)^{m+\alpha-1}\frac{(n_\epsi+1)^{\alpha-1} n_\epsi}{(1+\epsi n_\epsi)^3}\big(S_\epsi(\cdot,n_\epsi,c_\epsi)\nabla c_\epsi\cdot\nabla\varphi\big)\nonumber\\
&\qquad\qquad+\intoTomega(n_\epsi+1)^{m+2\alpha-1}(u_\epsi\cdot\nabla\varphi)\nonumber
\end{align}
holds for all $\epsi\in(0,1)$. Making use of \eqref{Sprop} we find that $\big|\frac{(n_\epsi+1)^{\alpha-1}n_\epsi}{(1+\epsi n_\epsi)^3}S_\epsi\big|\leq S_0$ for all $\epsi\in(0,1)$. Since moreover, $\frac{(n_\epsi+1)^{\alpha-1}n_\epsi}{(1+\epsi n_\epsi)^3}S_\epsi(\cdot,n_\epsi,c_\epsi)\to (n+1)^{\alpha-1}nS(\cdot,n,c)$ a.e. in $\Omega\times(0,\infty)$ as $\epsi\searrow0$ we find that 
\begin{align*}
\frac{(n_\epsi+1)^{\alpha-1}n_\epsi}{(1+\epsi n_\epsi)^3}S_\epsi(\cdot,n_\epsi,c_\epsi)\nabla c_\epsi\to (n+1)^{\alpha-1}nS(\cdot,n,c)\nabla c\quad\text{in }\LSp{2}{\Omega\times(0,T)}\quad\text{as }\epsi=\epsi_{j_k}\searrow0,
\end{align*}
in light of Lemma \ref{lem:strong-conv-nabc-l2} and \cite[Lemma 10.4]{win15_chemorot}. Combining this strong convergence with \eqref{eq:conv-n+1-m+a-1-ae} and \eqref{eq:conv-nab-n+1-m-1} entails that
\begin{align*}
\intoTomega(n_\epsi+1)^{m+\alpha-1}&\frac{(n_\epsi+1)^{\alpha-1} n_\epsi}{(1+\epsi n_\epsi)^3}\big(S_\epsi(\cdot,n_\epsi,c_\epsi)\nabla c_\epsi\cdot\nabla\varphi\big)\to\,
\intoTomega (n+1)^{m+2\alpha-2}n\big(S(\cdot,n,c)\nabla c\cdot\nabla\varphi\big)
\end{align*}
and
\begin{align*}
-
\intoTomega\frac{(n_\epsi+1)^{\alpha-1}n_\epsi}{(1+\epsi n_\epsi)^3}\big(\nabla(n_\epsi+1&)^{m+\alpha-1}\cdot S_\epsi(\cdot,n_\epsi,c_\epsi)\nabla c_\epsi\big)\varphi\\&\to-
\intoTomega(n+1)^{\alpha-1}n\big(\nabla (n+1)^{m+\alpha-1}\cdot S(\cdot,n,c)\nabla c\big)\varphi
\end{align*}
as $\epsi=\epsi_{j_k}\searrow0$, respectively. Moreover, relying on \eqref{eq:conv-n+1-m+2alpha-1}, \eqref{eq:conv-nab-n+1-m-1}, \eqref{eq:conv-n+1+eps-alpha-m-1} and \eqref{eq:conv-n+1-m+2alpha-1-u-l1} we obtain that
\begin{align*}
-\intoTomega(n_\epsi+1)^{m+2\alpha-1}\varphi_t&\to-\intoTomega (n+1)^{m+2\alpha-1}\varphi_t,\\
-\intoTomega(n_\epsi+1)^\alpha(n_\epsi+\epsi)^{m-1}\big(\nabla(n_\epsi+1)^{m+\alpha-1}\cdot\nabla\varphi\big)&\to-\intoTomega(n+1)^{\alpha}n^{m-1}\big(\nabla(n+1)^{m+\alpha-1}\cdot\nabla\varphi\big)\\
\intoTomega(n_\epsi+1)^{m+2\alpha-1}(u_\epsi\cdot\nabla\varphi)&\to\intoTomega n^{m+2\alpha-1}(u\cdot\nabla\varphi)
\end{align*}
as $\epsi=\epsi_{j_k}\searrow0$. Lastly, we depend on the lower semicontinuity of the norm in $\LSp{2}{\Omega\times(0,T)}$ with respect to weak convergence to conclude from \eqref{eq:conv-nab-n+1+eps-m-1} that
\begin{align*}
\liminf_{\epsi_{j_k}\searrow0}\intoTomega\big|(n_\epsi+1)^{\frac{m+2\alpha-3}{2}}(n_\epsi+\epsi)^{\frac{m-1}{2}}\nabla n_\epsi\big|^2\varphi\geq\intoTomega\big|(n+1)^\frac{m+2\alpha-3}{2}n^\frac{m-1}{2}\nabla n\big|^2\varphi.
\end{align*}
Uniting the statements above with \eqref{eq:sol-prop-n-proof-eq1} and the fact that $m+2\alpha<2$ entails that
\begin{align*}
&-\intinfomega(n+1)^{m+2\alpha-1}\varphi_t-\intomega(n_0+1)^{m+2\alpha-1}\varphi(\cdot,0)\\
\geq&\ m(m+2\alpha-1)(2-(m+2\alpha))\intinfomega\big|(n+1)^{\frac{m+2\alpha-3}{2}}n^{\frac{m-1}{2}}\nabla n\big|^2\varphi\\
&\quad-\frac{m(m+2\alpha-1)}{m+\alpha-1}\intinfomega(n+1)^\alpha n^{m-1}\big(\nabla(n+1)^{m+\alpha-1}\cdot\nabla\varphi\big)\\
&\qquad -\frac{(m+2\alpha-1)(2-(m+2\alpha))}{m+\alpha-1}\intinfomega (n+1)^{\alpha-1} n\big(\nabla(n+1)^{m+\alpha-1}\cdot S(\cdot,n,c)\nabla c\big)\varphi\\
&\quad\qquad+(m+2\alpha-1)\intinfomega(n+1)^{m+2\alpha-2}n\big(S(\cdot,n,c)\nabla c\cdot\nabla\varphi\big)+\intinfomega(n+1)^{m+2\alpha-1}(u\cdot\nabla\varphi),
\end{align*} 
where we used that $\varphi\equiv0$ in $\Omega\times(T,\infty)$. It can easily be checked, that this is an equivalent formulation of \eqref{eq:very_weak_supersol} for our choice $\Phi(s)\equiv(s+1)^{m+2\alpha-1}$, which thereby completes the proof. 
\end{bew}

The previous lemma at hand, we can conclude Theorem \ref{theo:2} in a straightforward manner.

\begin{proof}[\textbf{Proof of Theorem \ref{theo:2}}:]
The existence of a global very weak solution for $m\geq1$ and $\alpha\geq0$ satisfying $m+2\alpha>\frac53$ is already established in light of Theorem \ref{theo:1}, since any global weak solution is also a global very weak solution for the choice $\Phi(s)\equiv s$. So that clearly, we can restrict ourselves to verifying the $\Phi$--supersolution property of Definition \ref{def:weak_super_sol} for $m\geq1$, $\alpha\geq0$ satisfying $m+\alpha>\frac43$ and $m+2\alpha\leq\frac53$. In this case Lemma \ref{lem:sol-prop-n} is applicable and therefore, an evident combination of Lemmata \ref{lem:sol-prop-c-u} and \ref{lem:sol-prop-n} with the regularity information presented in Lemma \ref{lem:convergences} completes the proof.
\end{proof}

\section*{Acknowledgements}
The author acknowledges support of the {\em Deutsche Forschungsgemeinschaft} in the context of the project
  {\em Analysis of chemotactic cross-diffusion in complex frameworks}.

\footnotesize{
\setlength{\bibsep}{2pt plus 0.5ex}

}
\end{document}